\documentclass[11pt]{article}
\usepackage[english]{babel}
\usepackage[latin1]{inputenc}
\usepackage[T1]{fontenc}
\usepackage{amsthm,amsmath,amsbsy,amssymb,amsfonts,bbm}
\usepackage{graphicx}

\newcommand{\beq}{\begin{eqnarray}}
\newcommand{\eeq}{\end{eqnarray}}
\newcommand{\beqe}{\begin{eqnarray*}}
\newcommand{\eeqe}{\end{eqnarray*}}

\newcommand{\pa}[1]{\left({#1}\right)}
\newcommand{\cro}[1]{\left[{#1}\right]}

\newcommand{\ac}[1]{\left\{{#1}\right\}}
\newcommand{\LL}{\mathbb{L}}

\newcommand{\EE}{\mathbb{E}}
\newcommand{\PP}{\mathbb{P}}
\newcommand{\calS}{\mathcal{S}}

\newcommand{\N}{\mathbb{N}}

\newcommand{\1}{\mathbbm{1}}
\newtheorem{Th}{Theorem}

\newtheorem{Prop}{Proposition}
\newtheorem{Lemma}{Lemma}

\newtheorem{Cor}{Corollary}
\renewenvironment{proof}{\noindent{\bf Proof.}}{\hfill $\blacksquare$\par\noindent}

\newcommand{\norm}[1]{\ensuremath{\vert\!\vert #1 \vert\!\vert}}

\newcommand{\p}{\ensuremath{\varphi}}
\numberwithin{equation}{section}

\renewcommand{\sup}[1]{\mathop{\mathrm{sup}} \limits_{#1}}
\renewcommand{\inf}[1]{\mathop{\mathrm{inf}} \limits_{#1}}
\renewcommand{\max}[1]{\mathop{\mathrm{max}} \limits_{#1}}

\newcommand{\bJ}{{\bar{J}}}
\title{ Adaptive tests of homogeneity for a Poisson process}
\author{M. Fromont\thanks{ CREST-ENSAI / Equipe de Statistique de l'IRMAR, Université Rennes 2.} ,  B. Laurent\thanks{Institut de Mathématiques de Toulouse, INSA de Toulouse, Universit\'e de Toulouse.} , P. Reynaud-Bouret\thanks{CNRS ENS Paris et Universit\'e de Nice Sophia-Antipolis.}}

\begin{document}
\maketitle
\begin{abstract}
We propose to test the homogeneity of a Poisson process observed on a finite interval. In this framework, we first provide lower bounds for the uniform separation rates in $\mathbb{L}^2$ norm over classical Besov bodies and weak Besov bodies. Surprisingly, the obtained lower bounds over weak Besov bodies coincide with the minimax estimation rates over such classes. Then we construct non asymptotic and nonparametric testing procedures that are adaptive in the sense that they achieve, up to a possible logarithmic factor, the optimal uniform separation rates over various Besov bodies simultaneously. These procedures are based on model selection and thresholding methods. We finally complete our theoretical study with a Monte Carlo evaluation of the power of our tests under various alternatives.

\end{abstract}


\section{Introduction}

Poisson processes have been used for many years to model a great variety of situations: machine breakdowns, phone calls... Recently Poisson processes become popular for modeling occurrences of words or motifs on the DNA sequence (see  Robin, Rodolphe and Schbath \cite{RobinSchbath}). In this context, it is particularly important to be able to detect abnormal behaviors.

With such applications in mind, we consider in this paper the question of testing the homogeneity of a Poisson process $N$. Since we can only observe a finite number of points of the process, this question has a sense only on a finite interval. For the sake of simplicity, we assume that the Poisson process $N$ is observed on the fixed set $[0,1]$, and that it has an intensity $s$ with respect to some measure $\mu$ on $[0,1]$ with $d\mu(x)=L dx$.

Denoting by $\calS_0$ the set of constant
functions on $[0,1]$, our aim is consequently to test the null hypothesis $(H_0)$  "$s\in\calS_0$", against the alternative
$(H_1)$ "$s\not\in \calS_0$".

This problem of testing the homogeneity of a Poisson process has been widely investigated both from a theoretical and practical point of view (see Bain, Engelhardt, and Wright \cite{Bainetal} or Cohen and Sackrowitz \cite{Cohen} for a survey and Bhattacharjee, Deshpande, and Naik-Nimbalkar \cite{Bhatta} for a more recent work). In these papers, the alternative intensities are monotonous. Another related topic is the problem of testing the simple hypothesis that a stationary process is a Poisson process with a given intensity. We can cite for instance the papers by Fazli and Kutoyants \cite{FazliKut} where the alternative is also a Poisson process with a known intensity, Fazli \cite{Fazli} where the alternatives are Poisson processes with one-sided parametric intensities, or Dachian and Kutoyants \cite{DachianKut}, where the alternatives are self-exciting point processes. The paper by Ingster and Kutoyants \cite{IngsterKut} is the closest one to the present work. The alternatives considered by Ingster and Kutoyants are Poisson processes with nonparametric intensities in a Sobolev or Besov $\mathcal{B}_{2,q}^{(\delta)}(R)$ ball with $1\leq q <+\infty$ and known smoothness parameter $\delta$.

However, in some practical cases like the study of occurrences of words or motifs on a DNA sequence, such smooth alternatives cannot be considered. The intensity of the Poisson process in these cases may burst at a particular position of special interest for the biologist (see Gusto and Schbath \cite{Gusto} for more details). The question of testing the homogeneity of a Poisson process then becomes "how can we distinguish a Poisson process with constant intensity from a Poisson process whose intensity has some small localized spikes?". This question has already been partially considered in the seventies in a precursory work by Watson \cite{Watson}: he proposed a test based on the estimation of the Fourier coefficients of the intensity without evaluating the power of the resulting procedure.

In this paper, we focus on constructing adaptive testing procedures i.e. which do not use any prior information about the
smoothness of the intensity $s$, but which however have the best possible performances (in a minimax sense).

From a theoretical point of view, we evaluate the performances of the tests in terms of uniform separation rates
with respect to some prescribed distance $d$ over various classes of functions. Given $\beta\in]0,1[$, a class
of functions $\calS_1$, and a level $\alpha$ test $\Phi_{\alpha}$ with values in $\{0,1\}$ (rejecting $(H_0)$
when $\Phi_{\alpha}=1$), the uniform separation rate $\rho(\Phi_{\alpha},\calS_1,\beta)$ of $\Phi_{\alpha}$ over
the class $\calS_1$ is defined as the smallest positive number $\rho$ such that the test has an error of second
kind at most equal to $\beta$ for all alternatives $s$ in $\calS_1$ at an $\mathbb{L}^2$ distance $\rho$ from $\calS_0$. More
precisely, if $\PP_s$ denotes the distribution of the Poisson process $N$ with intensity $s$,

\begin{eqnarray}\label{defvitesse}
\rho(\Phi_\alpha,\calS_1,\beta)&=&\inf{}\left\{\rho>0, \sup{s\in\calS_1,
d(s,\calS_0)>\rho} \PP_s(\Phi_{\alpha}=0)\leq \beta\right\}\\
&=&\inf{}\left\{\rho>0,\inf{s\in\calS_1, d(s,\calS_0)> \rho} \PP_s(\Phi_{\alpha}=1)\geq 1-\beta\right\}.
\end{eqnarray}

In view of the practical situations of our interest, we study some classes of alternatives that can be very irregular, for instance that can have some localized spikes. We then consider some classical Besov bodies and also some spaces that can be viewed as weak versions of these classical Besov bodies and that are defined precisely in the following.  The interested reader may find in Rivoirard \cite{rivoirard2006} some illustrations of functions in weak Besov spaces and how the smoothness parameters of the functions govern the proportion and amplitude of their spikes.

As a first step, we evaluate the best possible value of the uniform separation rate over these spaces.
In other words, we give a lower bound for
\begin{eqnarray}\label{mrt}
\underline{\rho}(\calS_1,\alpha,\beta)= \inf{\Phi_{\alpha}}\rho(\Phi_{\alpha},\calS_1,\beta),
\end{eqnarray}
where the infimum is taken over all level $\alpha$ tests $\Phi_{\alpha}$, and where $\calS_1$ can be either a Besov body or a weak Besov body.
This quantity introduced by Baraud
\cite{yannick} as the $(\alpha,\beta)$-minimax rate of testing over $\calS_1$ or the minimax separation rate over
$\mathcal{S}_1$ is a stronger  version of the (asymptotic) minimax rate of testing usually considered.
The key reference for the computation of minimax rates of testing in various statistical models is
the series of papers due to Ingster \cite{Ing}. Concerning the Poisson model, Ingster and Kutoyants \cite{IngsterKut} give the minimax rate of testing for Sobolev or Besov $\mathcal{B}_{2,q}^{(\delta)}(R)$ balls with $1\leq q <+\infty$ and smoothness parameter $\delta>0$. They find that this rate of testing for the Sobolev or Besov norm or semi-norm is of order $L^{\frac{-2\delta}{4\delta+1}}$. Let us note that we find here lower bounds for the classical Besov bodies similar to Ingster and Kutoyants'ones. Furthermore, our lower bounds for the weak Besov bodies are larger than the ones for classical Besov bodies. Alternatives in weak Besov bodies are in fact so irregular that it is as difficult to detect them as to estimate them. The problem of estimation in weak Besov spaces is solved by using thresholding procedures: indeed the weak Besov spaces are closely related to the maxisets of those procedures (see  Kerkyacharian and Picard \cite{kerk} in the Gaussian framework and Reynaud-Bouret and Rivoirard \cite{PatouVincent} in a Poisson model). To our knowledge, no previous results of this kind exist for weak Besov bodies in testing problems, even in more classical statistical models, like the density model. Despite the similarity of both models, our lower bounds over weak Besov bodies cannot however be straightly transposed to the density model since our proofs heavily rely on the
 Poissonian independence properties.

As a second and main step, we construct non asymptotic level $\alpha$ tests which achieve, up to a possible logarithmic factor, the minimax separation rates over many Besov bodies and weak Besov bodies simultaneously, whereas using no prior information about the smoothness of the intensity $s$.
Our idea here is to combine some model selection methods that are effective for alternatives in classical Besov bodies and  a thresholding type approach, inspired by the thresholding rules used for adaptive estimation in weak Besov bodies. Key tools in the proofs of our results are exponential inequalities for U-statistics of order 2 due to Houdr\'e and Reynaud-Bouret \cite{ustats}.

Of course, both model selection and thresholding approaches have already been used to construct adaptive tests in various statistical models. One can cite among others the papers by Spokoiny (\cite{Spok1996} and \cite{Spok1998}) in Gaussian white noise models or by Baraud, Huet and Laurent \cite{BHL} in a Gaussian regression framework. These papers propose adaptive tests which combine methods closely related to both model selection and thresholding ones. As for the density framework, adaptive tests were proposed by Ingster \cite{Ing2000} or Fromont and Laurent \cite{Beamag}, using model selection type methods and by Butucea and Tribouley \cite{Karine} using thresholding type methods.

The present work is organized as follows. In Section \ref{sectionmino}, we provide lower bounds  for the uniform separation rates over various Besov bodies. Our testing procedures are defined in Section \ref{sectionproc}, and their uniform separation rates over Besov bodies are established in Section \ref{sectionvitesses}. We carry out a simulation study in Section \ref{sectionsimu} to illustrate these theoretical results, and the proofs are postponed to the last section.

\section{Lower bounds for the minimax separation rates over Besov bodies}\label{sectionmino}

We consider the Poisson process $N$ with intensity $s$ with respect to some measure $\mu$ on $[0,1]$, with
$d\mu(x)=Ldx$. In the following, we assume that $s$ belongs to $\mathbb{L}^2([0,1])$, and $\langle .,.\rangle$, $\|.\|$ and $d$ respectively denote the scalar product
$$\langle f,g\rangle=\int_{[0,1]} f(x)g(x) dx,$$
the $\mathbb{L}^2-$norm
$$\|f\|^2=\int_{[0,1]}f^2(x)dx,$$  and the associated distance.

Let us denote the Haar basis of $\mathbb{L}^2([0,1])$ by $\{\phi_0,\phi_{(j,k)},j\in\N,k\in\{0,\ldots,2^j-1\}\}$ with
$$\phi_0(x)=\1_{[0,1]}(x),$$
and
\begin{equation}\label{Haar}
\phi_{(j,k)}(x)=2^{j/2}\psi(2^jx-k),
\end{equation}
where $\psi(x)=\1_{[0,1/2[}(x)-\1_{[1/2,1[}(x)$.

We set $\alpha_0=\langle s,\phi_0\rangle$ and for every $j\in\N,k\in\{0,\ldots,2^j-1\}$, $\alpha_{(j,k)}=\langle s,\phi_{(j,k)}\rangle$.

We can now introduce the Besov bodies defined for $\delta>0$, $R>0$ by
\begin{multline}\label{defbesovfort}
 \mathcal{B}_{2,\infty}^\delta(R)=\Bigg\{s\geq 0,\
s\in\LL^2([0,1]),
 s=\alpha_0\phi_0+\sum_{j\in\N}\sum_{k=0}^{2^j-1} \alpha_{(j,k)} \phi_{(j,k)},\\ \forall j
\in\N, \sum_{k=0}^{2^j-1} \alpha_{(j,k)}^2\leq R^2 2^{-2j\delta}\Bigg\},
\end{multline}

and more generally for $p\geq 1$, $R>0$ and $\delta>\max~(0, 1/p-1/2)$,

\begin{multline}
 \mathcal{B}_{p,\infty}^\delta(R)=\Bigg\{s\geq 0,\
s\in\LL^2([0,1]),
 s=\alpha_0\phi_0+\sum_{j\in\N}\sum_{k=0}^{2^j-1} \alpha_{(j,k)} \phi_{(j,k)},\\ \forall j
\in\N, \sum_{k=0}^{2^j-1}|\alpha_{(j,k)}|^p\leq R^p 2^{-pj\left(\delta+\frac{1}{2}-\frac{1}{p}\right)}\Bigg\}.
\end{multline}

As in Reynaud-Bouret and Rivoirard \cite{PatouVincent}, we also introduce some weaker versions of the above Besov bodies given for $\gamma >0$ and $R'>0$ by
\begin{multline}\label{defbesovfaible}
W_{\gamma}(R')=\Bigg\{s \geq 0, s \in \mathbb{L}^2([0,1]), s=\alpha_0\phi_0+\sum_{j \in \mathbb{N}}\sum_{k=0}^{2^j-1}\alpha_{(j,k)}\phi_{(j,k)},\\
 \forall t>0, \sum_{j
\in\N} \sum_{k=0}^{2^j-1} \alpha_{(j,k)}^2\1_{\alpha_{(j,k)}^2\leq t}\leq R'^2  t^{\frac{2\gamma}{1+2\gamma}}\Bigg\}.
\end{multline}

Fixing some levels of error $\alpha$ and $\beta$ in $]0,1[$, and denoting by $\mathbb{L}^\infty(R'')$ the set of functions bounded by $R''$, our purpose in this section is to find sharp lower bounds for $\underline{\rho}(\mathcal{B}_{2,\infty}^\delta(R)\cap W_{\gamma}(R')\cap \mathbb{L}^{\infty}(R''),\alpha,\beta)$,
where $\underline{\rho}$ is defined by (\ref{mrt}).

Starting from a general idea developed by Ingster \cite{Ing}, we obtain the following result.

\begin{Th}\label{minorations}
Assume that $R>0$, $R'>0$, and  $R''\geq 2$, and fix some levels $\alpha$ and $\beta$ in $]0,1[$ such that $\alpha+\beta\leq 0.59$.

$(i)$ If $\gamma > \max{}(2\delta,1/2)$,  then
$$\liminf_{L\to +\infty}\left(\frac{L}{\ln L}\right)^{\frac{\gamma}{1+2\gamma}}\underline{\rho}( \mathcal{B}_{2,\infty}^{\delta}(R)\cap W_{\gamma}(R')\cap \mathbb{L}^{\infty}(R''),\alpha,\beta)>0 .$$
$(ii)$ If $\delta \geq \max{}(\gamma/2,\gamma/(1+2\gamma))$, then
$$\liminf_{L\to +\infty}L^{\frac{2\delta}{1+4 \delta}}\underline{\rho}( \mathcal{B}_{2,\infty}^{\delta}(R)\cap W_{\gamma}(R')\cap \mathbb{L}^{\infty}(R''),\alpha,\beta)>0.$$
$(iii)$ If $\delta<\gamma/(1+2\gamma)$ and $\gamma\leq 1/2$, then
$$\liminf_{L\to +\infty}L^{\frac{1}{2}}\underline{\rho}( \mathcal{B}_{2,\infty}^{\delta}(R)\cap W_{\gamma}(R')\cap \mathbb{L}^{\infty}(R''),\alpha,\beta)>0.$$
\end{Th}

\emph{Comments.}

\begin{enumerate}
\item  For the whole set of parameters $(\delta,\gamma)$ such that $\delta \geq \gamma/(1+2\gamma)$, we prove in Section \ref{sectionvitesses} that these lower bounds are actually sharp.
\item We have in case $(ii)$ lower bounds  which coincide with the minimax rates of testing obtained by Ingster and Kutoyants \cite{IngsterKut} when testing that a periodic Poisson process has a given intensity in the Besov spaces $ \mathcal{B}_{2,\infty}^{(\delta)}(R)$. We know (see Ingster \cite{Ing2000} or Fromont and Laurent \cite{Beamag} for instance) that such rates can be achieved by some multiple testing procedure based on model selection type methods. This is the principle of our first procedure described in Section \ref{procedureselmod}.
\item  We notice that the lower bounds obtained in case $(i)$ are equal to the minimax estimation rates on the maxisets of the thresholding estimation procedure, namely $\mathcal{B}_{2,\infty}^{\gamma/(1+2\gamma)}(R)\cap W_{\gamma}(R')$ (see Kerkyacharian and Picard \cite{kerk}, Rivoirard \cite{rivoirard2006}, or Reynaud-Bouret and Rivoirard \cite{PatouVincent}  for more details). This means that it is as difficult to test as to estimate over such classes of functions, phenomenon  which is quite unusual. Since the minimax estimation rates on these classes are achieved  by thresholding rules, it will be natural to construct a testing procedure based on thresholding methods: this is the idea that originated our second procedure described in Section \ref{procedureseuillage}.
\end{enumerate}

\section{Two tests of homogeneity}\label{sectionproc}

Let us recall that $\calS_0$ denotes the set of constant functions on $[0,1]$ and that we assume that $s$ belongs to $\mathbb{L}^2([0,1])$.

In this section, we construct level $\alpha$ tests of the null
hypothesis $(H_0)$ "$s\in\calS_0$", against the alternative $(H_1)$  "$s\not\in \calS_0$", from the observation
of the Poisson process $N$, or the points $\{X_l,\  l=1,\ldots, N_L\}$ of the Poisson process.

We introduce two testing procedures that come from two different statistical approaches. The first one
originates in general model selection methods, while the second one is closer to the thresholding type methods.

In order to understand the global ideas of these procedures, let us notice that the
squared $\mathbb{L}^2-$distance $d^2(s,\calS_0)$ between $s$ and the set of constant functions $\calS_0$ can be rewritten as \begin{eqnarray*} d^2(s,\calS_0)&=&\int_{[0,1]}\pa{s(x)- \int_{[0,1]} s(y)dy}^2 dx,\nonumber\\
&=& \|s\|^2-\alpha_0^2,\nonumber\\
&=&\sum_{\lambda\in\Lambda_{\infty}} \alpha_{\lambda}^2,\end{eqnarray*}
where $\alpha_0=\langle s,\phi_0 \rangle$, and for all $\lambda\in \Lambda_{\infty}=\{(j,k),\ j\in\N,\ k\in\{0,\ldots,2^j-1\}\}$, $\alpha_{\lambda}=\langle
 s,\phi_{\lambda} \rangle $.

For every $\lambda \in
  \Lambda_{\infty}$, $\alpha_\lambda$ can be estimated by
$$ \widehat{\alpha_\lambda}=\frac{1}{L} \int_{[0,1]} \phi_\lambda(x) dN_x,$$
which is also equal to
$$ \widehat{\alpha_\lambda}=\frac{1}{L} \sum_{l=1}^{N_L} \phi_\lambda(X_l).$$
From this variable, we deduce an unbiased estimator of $\alpha_\lambda^2$ given by :
\beq\label{estcoeff}
T_\lambda=\widehat{\alpha_\lambda}^2-\frac{1}{L^2} \int_{[0,1]}\phi_\lambda^2(x) dN_x=\frac{1}{L^2} \sum_{l \neq
l'=1}^{N_L} \phi_\lambda(X_l)\phi_\lambda(X_{l'}).\eeq

Our first approach will consist in constructing estimators of $d^2(s,\calS_0)=\sum_{\lambda\in\Lambda_{\infty}} \alpha_{\lambda}^2$ based on a combination of the $T_{\lambda}$'s, and in rejecting the null hypothesis when one of these estimators is too large. This was already the spirit of  Watson's procedure (see \cite{Watson}).
Our second approach is related to the test considered in Baraud et al. \cite{BHL} to detect local alternatives.  It will consist in considering a set of $T_\lambda$'s and rejecting the null hypothesis directly when one of the $T_{\lambda}$'s is too large. Let us now precisely define both procedures.

\subsection{A first procedure based on model selection}\label{procedureselmod}

Assuming that $s\in\LL^2([0,1])$, a natural idea is to construct a testing procedure from an estimation of the
squared $\mathbb{L}^2-$distance $d^2(s,\calS_0)$.

In order to estimate this functional of $s$, following the ideas of Laurent \cite{Bea2003} and Fromont and Laurent \cite{Beamag}, we
 introduce embedded finite dimensional linear subspaces of $\LL^2([0,1])$. We choose here to consider for $J\geq 1$ the subspaces $S_J$ generated by the subsets $\{\phi_0,\phi_\lambda, \lambda \in \Lambda_J\}$ of the Haar basis defined by (\ref{Haar}), with $\Lambda_J=\{(j,k),\ j\in\{0,\ldots, J-1\},\ k\in\{0,\ldots,2^j-1\}\}$. Each subspace $S_J$ is called a \emph{model}. We denote by $D_J=2^J$ the dimension of $S_J$, and by $s_J$ the orthogonal projection of $s$ onto the model $S_J$.

Focusing on one model $S_J$, we estimate
$d^2(s,\calS_0)=\|s\|^2-\alpha_0^2$ by the unbiased estimator of $\|s_J\|^2-\alpha_0^2=\sum_{\lambda\in \Lambda_J} \alpha_\lambda^2$ given by
\beq\label{estglob} T_J'=\sum_{\lambda \in
\Lambda_J}T_\lambda,\eeq
with $T_{\lambda}$ defined by (\ref{estcoeff}).
The estimator $T_J'$ obviously depends
on the choice of the model $S_J$.

Since we do not want to choose a priori such a model, we consider a collection
of models $\{S_J, J\in \mathcal{J}\}$ where $\mathcal{J}$ is a finite subset of $\N^*$, and the corresponding collection of estimators $\{T_J', J\in \mathcal{J}\}$.

The procedure that we introduce here then consists in rejecting $(H_0)$ "$s\in \calS_0$" when there
exists $J$ in $\mathcal{J}$ such that the estimator $T_J'$ given by (\ref{estglob}) is too large.

At this point there are several ways to decide when $T_J'$ is too large.

In all cases, we use the well-known argument that, conditionally on the event "the number of points $N_L$ falling into $[0,1]$
is $n$", the points of the process obeys the same law as a $n$-sample $(\tilde{X}_1, \ldots, \tilde{X}_n)$ with
common density $s/ \int_{[0,1]}s(x) dx$. It follows that for all $n \in \N$,
$$ \PP_s\pa{T_J'>q'|N_L=n}=\PP \pa{ \frac{1}{L^2} \sum_{\lambda \in
\Lambda_J} \sum_{l\neq l'=1}^n \phi_\lambda(\tilde{X}_l)\phi_\lambda(\tilde{X}_{l'})>q'}.$$

Under the null hypothesis, the intensity $s$ is constant on $[0,1]$, and the $\tilde{X}_l$'s are i.i.d., with
uniform distribution on $[0,1]$. This distribution is free from the parameter $s$. As a consequence, for every
$u\in]0,1[$, we can introduce and estimate by Monte Carlo experiments the $(1-u)$ quantile of the distribution
of $T_J'|N_L=n$ under the null hypothesis, that we denote by $q_J'^{(n)}(u)$.

We now consider the test statistics:

 \beq \label{stattest1} \mathcal{T}_{\alpha}^{(1)}=\sup{J\in\mathcal{J}}\left(T_J'-q_{J}'^{(N_L)}(u_{J,\alpha}'^{(N_L)})\right), \eeq
with $u_{J,\alpha}'^{(N_L)}$ to be
correctly chosen.

Finally, we define the corresponding test function:

\beq \label{fonctiontest1} \Phi_{\alpha}^{(1)}=\1_{\mathcal{T}_{\alpha}^{(1)}>0}. \eeq

And our first test consists in rejecting the null
hypothesis $(H_0)$ when $\Phi_{\alpha}^{(1)}=1$.\\

Let us see how we can choose $u_{J,\alpha}'^{(N_L)}$ so that our test has a level $\alpha$.

An obvious possibility  is to set
$$ u_{J,\alpha}'^{(n)}=\frac{\alpha}{|\mathcal{J}|}\textrm{  for every }J \textrm{ in } \mathcal{J} \textrm{ and } n \textrm{ in } \N.$$
This choice corresponds to a Bonferroni procedure and $\Phi_{\alpha}^{(1)}$ actually defines a level $\alpha$
test. Indeed, for $s\in\calS_0$,
\begin{eqnarray*}
&&\PP_s\pa{ \sup{J\in\mathcal{J}}\left(T_J'-q_J'^{(N_L)}(u_{J,\alpha}'^{(N_L)})\right)>0 }\\
&&=\sum_{n\in \N} \PP_s\pa{\sup{J\in\mathcal{J}}\left(T_J'-q_J'^{(n)}\left(\frac{\alpha}{|\mathcal{J}|}\right)\right)>0\Bigg|N_L=n}\PP_s(N_L=n)\\
&&\leq\sum_{n\in\N}\sum_{J\in\mathcal{J}}\frac{\alpha}{|\mathcal{J}|}\PP_s(N_L=n)\\
&&\leq \alpha.
\end{eqnarray*}

Our choice for $u_{J,\alpha}'^{(n)}$, inspired by Fromont and Laurent \cite{Beamag}, leads to a less conservative procedure. It consists in setting \beq
\label{choixpoids} u_{J,\alpha}'^{(n)}= e^{-W_J}\sup{}\ac{u \in ]0,1[, \sup{s\in\calS_0}\PP_{s}\pa{\sup{J\in\mathcal{J}}\left(
T_J'-q_J'^{(n)}(ue^{-W_J})\right)>0\Bigg|N_L=n} \leq \alpha }, \eeq where $\ac{W_J,J\in\mathcal{J}}$ is a collection of
positive weights such that
$$\sum_{J\in\mathcal{J}} e^{-W_J}\leq 1.$$

For the same reason, we still obtain a level $\alpha$ test and by definition, $u_{J,\alpha}'^{(n)}\geq
\alpha e^{-W_J}$ for every $n$ in $\N$.\\

\subsection{A second procedure based on a thresholding approach}\label{procedureseuillage}

Let us recall here that the
squared $\mathbb{L}^2-$distance $d^2(s,\calS_0)$ between $s$ and the set of constant functions $\calS_0$ is equal to $\sum_{\lambda\in\Lambda_{\infty}} \alpha_{\lambda}^2$ and that $T_\lambda$ defined by (\ref{estcoeff}) is an unbiased  estimator of $\alpha_{\lambda}^2$. Based on general thresholding ideas, our second procedure consists in fixing some $\bJ\geq 1$ and rejecting the null hypothesis $(H_0)$ when there exists $\lambda$ in $\Lambda_\bJ$ such that $T_\lambda$ is too large.

Let us now see what we mean by "$T_\lambda$ is too large". We can still
use the fact that
$$ \PP_s\pa{T_\lambda>q|N_L=n}=\PP \pa{ \frac{1}{L^2} \sum_{l \neq l'=1}^{n}
\phi_\lambda(\tilde{X}_l)\phi_\lambda(\tilde{X}_{l'})>q},$$ and that under the null hypothesis, the
$\tilde{X}_l$'s are i.i.d., with uniform distribution on $[0,1]$. We therefore introduce and estimate by Monte
Carlo experiments the $(1-u)$ quantile of the distribution of $T_\lambda|N_L=n$ under the null hypothesis, that
we denote by $q_\lambda^{(n)}(u)$.  Notice that for $\lambda=(j,k)\in \Lambda_\bJ$, $q_\lambda^{(n)}(u)$ does not
depend on $k$.

We set \beq \label{stattest2} \mathcal{T}_{\alpha}^{(2)}=\sup{\lambda\in\Lambda_\bJ}\left(T_\lambda-q_{\lambda}^{(N_L)}(u_{\lambda,\alpha}^{(N_L)})\right), \eeq
with
 $u_{\lambda,\alpha}^{(N_L)}$ to be correctly chosen.

We also define
\beq \label{fonctiontest2} \Phi_{\alpha}^{(2)}=\1_{\mathcal{T}_{\alpha}^{(2)}>0}. \eeq

 Our test consists in rejecting the null
hypothesis $(H_0)$ when $\Phi_{\alpha}^{(2)}=1$.\\

Let us now see how we choose $u_{\lambda,\alpha}^{(n)}$. An obvious choice corresponding to the Bonferroni procedure
would be
$$ u_{(j,k),\alpha}^{(n)}=\frac{\alpha}{2^j \bJ}\textrm{  for every }(j,k) \textrm{ in } \Lambda_\bJ \textrm{ and } n \textrm{ in } \N.$$
To obtain a less conservative procedure, we prefer setting
\beq \label{choixseuillage}
u_{\lambda,\alpha}^{(n)} = \frac{u_{\alpha}^{(n)}}{2^j \bJ}\textrm{  for every }\lambda=(j,k) \textrm{ in } \Lambda_\bJ \textrm{ and } n \textrm{ in } \N,\eeq with
 \beq
\label{ualpha}
u_{\alpha}^{(n)} = \sup{}  \left\{ u \in ]0,1[, \sup{s\in\calS_0}\PP_{s} \pa{ \sup{(j,k)\in\Lambda_\bJ} \left(T_{(j,k)}-q_{(j,k)}^{(n)} \pa{\frac{u}{2^j \bJ}} \right)>0 \Bigg|N_L=n }\leq \alpha \right\} .
\eeq
When $s\in\calS_0$,
 \begin{eqnarray*}
&&\PP_s\pa{ \sup{\lambda\in \Lambda_\bJ}\left(T_\lambda-q_\lambda^{(N_L)}\left(u_{\lambda,\alpha}^{(N_L)}\right)\right)>0}\\
&&= \PP_s\pa{ \sup{(j,k)\in\Lambda_\bJ}\left(T_{(j,k)}-q_{(j,k)}^{(N_L)}\left(\frac{u_{\alpha}^{(N_L)}}{2^{j}\bJ}\right)\right)>0}\\
&&= \sum_{n\in\N}\PP_s\pa{ \sup{(j,k)\in\Lambda_\bJ}\left(T_{(j,k)}-q_{(j,k)}^{(n)}\left(\frac{u_{\alpha}^{(n)}}{2^{j}\bJ}\right)\right)>0\Bigg|N_L=n}\PP(N_L=n)\\
&&\leq \alpha,
 \end{eqnarray*}
which means that $\Phi_{\alpha}^{(2)}$ defines a level $\alpha$ test.\\
Note that $u_{\alpha}^{(n)} \geq \alpha$.\\

\emph{Comments.}

\begin{enumerate}

\item Though the two testing procedures defined by (\ref{fonctiontest1}) and (\ref{fonctiontest2}) are very different by their spirit, they can formally be written in a common way. For any subset $\Lambda$ of $\Lambda_{\infty}$, we denote by $S_{\Lambda}$ the subspace generated by $\{\phi_{0},\phi_{\lambda},\lambda \in \Lambda\}$, by $s_\Lambda$ the orthogonal projection of $s$ onto $S_\Lambda$, and we introduce the unbiased estimator $T_\Lambda''=\sum_{\lambda\in \Lambda} T_\lambda$  of $\|s_\Lambda\|^2-\alpha_0^2=\sum_{\lambda\in \Lambda} \alpha_{\lambda}^2$. Then our test functions can be written as
\beq \label{fonctiontestgen}
\Phi_{\alpha}=\1_{\mathcal{T}_{\alpha}>0},
\eeq
where $$\mathcal{T}_{\alpha}=\sup{\Lambda\in\mathcal{C}}\left(T_\Lambda''-{t_{\Lambda,\alpha}''}^{(N_L)}\right),$$
and $\mathcal{C}$ is a finite collection of subsets of $\Lambda_{\infty}$.

Noticing that $T_J'=T_{\Lambda_J}''$, we can easily see that our first test amounts in taking a collection $\mathcal{C}$ equal to $\{\Lambda_J,\ J\in\mathcal{J}\}$, and ${t_{\Lambda_J,\alpha}''}^{(N_L)}=q_{J}'^{(N_L)}(u_{J,\alpha}'^{(N_L)})$.
Furthermore, our second test amounts in taking a collection $\mathcal{C}$ composed of all subsets of $\Lambda_\bJ$, and for $\Lambda\subset\Lambda_\bJ$, ${t_{\Lambda,\alpha}''}^{(N_L)}=\sum_{\lambda\in \Lambda} q_{\lambda}^{(N_L)}(u_{\lambda,\alpha}^{(N_L)})$.
Indeed, there exists a subset $\Lambda$ of $\Lambda_\bJ$ such that $\sum_{\lambda\in\Lambda}T_{\lambda}>\sum_{\lambda\in\Lambda}q_{\lambda}^{(N_L)}(u_{\lambda,\alpha}^{(N_L)})$ if and only if there exists $\lambda$ in $\Lambda_\bJ$ such that $T_{\lambda}>q_{\lambda}^{(N_L)}(u_{\lambda,\alpha}^{(N_L)})$.

Such a common expression will be particularly useful to derive the properties of the tests.

It also allows us to see our tests as multiple testing procedures. Indeed, we can consider that for each $\Lambda$ in $\mathcal{C}$, we construct a test rejecting the null hypothesis when $T_\Lambda''-{t_{\Lambda,\alpha}''}^{(N_L)}>0$. We thus obtain a collection of tests and we finally decide to reject the null hypothesis when it is rejected for at least one of the tests of the collection.

\item Both procedures have a specific interest to prove the optimality of the lower bounds obtained in Theorem \ref{minorations}. We will actually prove in the next section that the first one achieves the lower bounds obtained in case $(ii)$ of Theorem \ref{minorations} (up to a possible logarithmic factor) whereas the second one achieves the lower bounds obtained in case $(i)$ of Theorem \ref{minorations}. However, if we want a procedure that achieves the lower bounds of cases $(i)$ and $(ii)$ simultaneously, we will have to consider the test which consists in mixing the two procedures. In this case, we reject the
null hypothesis $(H_0)$ when $\sup{}\{\Phi_{\alpha/2}^{(1)},\Phi_{\alpha/2}^{(2)}\}=1$.

\end{enumerate}

\section{Uniform separation rates}\label{sectionvitesses}

In this section, we evaluate the performances of our new testing procedures from a theoretical point of view. More precisely, we prove that our procedures are optimal in the sense that their uniform separation rates over Besov bodies  are of the same order as the lower bounds for $\underline{\rho}$ obtained in Section \ref{sectionmino}. These results justify the construction of our procedures as well as they provide the upper bounds needed for the exact evaluation of the minimax separation rates over weak and classical Besov bodies in the Poisson framework.

\smallskip

In the following, the expression $C(\alpha, \beta,R, R',R'',\delta,\gamma,...)$ or $C_k(\alpha,\beta,R, R',R'',\delta,\gamma...)$ is used to denote some constant which only depends on the parameters $\alpha, \beta, R, R',R'',\delta,\gamma,...$, and which may vary from line to line.

\subsection{Uniform separation rates of the first procedure}

\subsubsection{The error of second kind}

The aim of the following theorem is to give a condition on the alternative so that our first level $\alpha$ test
has a prescribed error of second kind.

\begin{Th}\label{puissancetest1}
 Assume that $s\in\LL^{\infty}([0,1])$, and that $L\geq 1$. Fix some levels $\alpha$ and $\beta$
in $]0,1[$, and let $\Phi_{\alpha}^{(1)}$ be the test function defined by
(\ref{fonctiontest1}). There exist some positive constants
$C_1(\norm{s}_\infty,\beta)$, $C_2(\beta)$, $C_3(\alpha,\beta)$, and
$C_4(\alpha)$ such that when $s$ satisfies

\begin{multline}\label{oracle}
d^2(s,\calS_0) > \inf{J\in\mathcal{J}}\Bigg\{ \norm{s-s_J}^2
+C_1(\norm{s}_\infty,\beta)\frac{\sqrt{D_J}}{L}+C_2(\beta)\frac{D_J}{L^2}\\
+C_3(\alpha,\beta)\int_{[0,1]}
s(x)dx\pa{\frac{\sqrt{D_JW_J}}{L}+\frac{W_J}{L}} + C_4(\alpha)\frac{D_J
W_J^2}{L^2}\Bigg\}\end{multline}

then
$$ \PP_s\pa{\Phi_{\alpha}^{(1)}=0}\leq \beta. $$

\end{Th}

\emph{Comment.}  Considering here a multiple testing procedure instead of a simple one allows to obtain in the right hand side of the inequality $(\ref{oracle})$ an infimum over all $J$ in $\mathcal{J}$ at the only price of introducing some terms in $W_J$. These last terms will appear in the following uniform separation rates over classical Besov bodies as a $\ln \ln L$ factor, which is now known to be the price to pay for adaptivity in some classical statistical models. As a consequence, our multiple testing procedure is proved to be adaptive in Proposition \ref{vitesses1} over classical Besov bodies, which would not occur with a simple testing procedure.

\subsubsection{Uniform separation rates over Besov bodies}

In this section, we evaluate the uniform separation rates $\rho(\Phi_\alpha^{(1)},\mathcal{B}_{2,\infty}^\delta(R)\cap\LL^{\infty}(R''),\beta)$
where $\rho$ is defined by (\ref{defvitesse}), and $\mathcal{B}_{2,\infty}^\delta(R)$ is any Besov body defined by $(\ref{defbesovfort})$.

 Let us first notice that the functions of $\mathcal{B}_{2,\infty}^\delta(R)$ are well approximated by their projections onto subspaces of the collection
$\{S_J,\ J\in\mathcal{J}\}$ considered in our first procedure, in the sense that if
$s\in\mathcal{B}_{2,\infty}^\delta(R)$, then
$$\|s-s_J\|^2\leq c(\delta) R^2 D_J^{-2\delta}.$$

As a consequence we can use Theorem \ref{puissancetest1} to obtain upper bounds for the uniform separation rates
of our test. \\
We denote by $\left\lfloor x \right \rfloor$ the integer part of $x$.

\begin{Prop}\label{vitesses1}

Assume that $\ln \ln L\geq 1$. Given some levels $\alpha$ and $\beta$ in $]0,1[$, let
$\Phi_\alpha^{(1)}$ defined by (\ref{fonctiontest1}) with $\mathcal{J}=\{1,\ldots,\left\lfloor\log_2(L^2/(\ln\ln L)^3)\right \rfloor\}$ and $W_J=\ln|\mathcal{J}|$ for every $J$ in $\mathcal{J}$.

For every
$\delta>0$, there exists some positive constant $C(\alpha,\beta,R'',\delta)$ such that when $s$ belongs to $\mathcal{B}_{2,\infty}^\delta(R)\cap \LL^\infty(R'')$ and satisfies
$$d^2(s,\calS_0)>
C(\alpha,\beta,R'',\delta)\left(R^{\frac{2}{4\delta +1}} \left( \frac{\sqrt{\ln\ln
L}}{L}\right)^{\frac{4\delta}{4\delta +1}}+R^2\left(\frac{(\ln \ln L)^3}{L^2}\right)^{2\delta}+\frac{\ln\ln L}{L}\right),$$
then
$$\PP_s\pa{\Phi_{\alpha}^{(1)}=0}\leq \beta. $$
In particular, there exist some positive constants $L_0(\delta)$ and $C(\alpha,\beta,R,R'',\delta)$ such that if $L>L_0(\delta)$, then
$$\rho(\Phi_\alpha^{(1)},\mathcal{B}_{2,\infty}^\delta(R)\cap \LL^\infty(R''),\beta)\leq
C(\alpha,\beta,R,R'',\delta)\left( \frac{\sqrt{\ln\ln
L}}{L}\right)^{\frac{2\delta}{4\delta +1}}.$$

\end{Prop}

\emph{Comments.}
\begin{enumerate}
\item Our first testing procedure is therefore adaptive: indeed, for large $L$, it achieves the lower bounds for the minimax separation
rates over all the spaces $\mathcal{B}_{2,\infty}^\delta(R)\cap W_\gamma(R')\cap \LL^\infty(R'')$ with $\delta \geq \max{}(\gamma/{2},\gamma/(1+2\gamma))$ simultaneously up to a possible $\ln\ln L$ factor (see Theorem~\ref{minorations}).\\
However it does not achieve the optimal separation rates obtained in the case where $\gamma > \max{}(2\delta,1/2)$. In this range of parameters, the regularity in $\gamma$ is higher than the regularity in $\delta$, meaning that the weak Besov body governs the separation rate.  That is the reason why we introduced the thresholding type procedure.

\item The upper bound for the uniform separation rate obtained here is exactly of the same order as the (asymptotic)
adaptive minimax rate of testing obtained by Ingster \cite{Ing2000} in the density model, replacing the parameter $L$ of the Poisson model  by the number $n$ of observations in the density model. In particular the $\ln \ln L$ factor is proved to be necessary in the density model for adaptive procedures.

\item It is easy to see that $\mathcal{B}_{p,\infty}^\delta(R)\subset \mathcal{B}_{2,\infty}^{\delta}(R)$ when $p>2$. So this result directly leads to upper bounds for the uniform separation rates $\rho(\Phi_\alpha^{(1)},\mathcal{B}_{p,\infty}^\delta(R)\cap \LL^\infty(R''),\beta)$ of the test over the Besov bodies $\mathcal{B}_{p,\infty}^\delta(R)$ when $p>2$. These rates, obtained in the Poisson framework, correspond to the ones in some Gaussian models (see Spokoiny \cite{Spok1996} for instance) or in the density model (see Ingster \cite{Ing2000}).

\item Note that one could also consider some tests based on the Fourier basis as well as the Haar basis, as Fromont and Laurent \cite{Beamag} did in the density model. The theoretical results would remain unchanged, and the practical performances of the procedure would be better when considering smooth alternatives (see Fromont and Laurent  \cite{Beamag} for more details and Section \ref{sectionsimu}). We have only considered here tests based on the Haar basis for the sake of simplicity.

\end{enumerate}

\subsection{Uniform separation rates of the second procedure}

\subsubsection{The error of second kind}

From the common expression (\ref{fonctiontestgen}) of the test function for the two procedures, we obtain here a result similar to Theorem \ref{puissancetest1} for the error of second kind of our second test.

\begin{Th}\label{puissancetest2}
Assume that
$s\in\LL^{\infty}([0,1])$, and that $L\geq 1$. Fix some levels $\alpha$ and $\beta$ in $]0,1[$, and let $\Phi_{\alpha}^{(2)}$ be the test function defined by (\ref{fonctiontest2}). Recall that for any subset $\Lambda$ of $\Lambda_{\infty}$, $S_{\Lambda}$ and $s_{\Lambda}$ respectively denote the subspace generated by $\{\phi_{0},\phi_{\lambda},\lambda \in \Lambda\}$ and the orthogonal projection of $s$ onto $S_\Lambda$. Denoting by $D_\Lambda$ the dimension of $S_\Lambda$, there exist some
positive constants $C_1(\norm{s}_\infty,\beta)$, $C_2(\beta)$,
$C_3(\alpha,\beta)$, and $C_4(\alpha)$ such that when $s$ satisfies
\begin{multline}
d^2(s,\calS_0) > \inf{\Lambda\subset \Lambda_\bJ}\Bigg\{ \norm{s-s_\Lambda}^2
+C_1(\norm{s}_\infty,\beta)\pa{\frac{\sqrt{D_\Lambda}}{L}+\frac{2^{\bJ/2}}{L^{3/2}}}+C_2(\beta)\frac{2^\bJ}{L^2}\\
+C_3(\alpha,\beta)\int_{[0,1]} s(x)dx\frac{D_\Lambda\ln(2^{\bJ}\bJ)}{L} + C_4(\alpha)\frac{D_\Lambda 2^\bJ \ln^2(2^\bJ \bJ)}{L^2}\Bigg\}\end{multline}

then
$$ \PP_s\pa{\Phi_{\alpha}^{(2)}=0}\leq \beta. $$

\end{Th}

\subsubsection{Uniform separation rates over Besov bodies}

\begin{Prop}\label{vitesses2}

Assume that $\ln L\geq 1$. Given some levels $\alpha$ and $\beta$ in $]0,1[$, let  $\Phi_{\alpha}^{(2)}$ be the test defined by (\ref{fonctiontest2}) with $\bJ=\left\lfloor\log_2(L/\ln L)\right\rfloor$.

For every $\delta>0$ and $\gamma>0$, there exists some positive constant $C(\alpha,\beta,R'',\delta,\gamma)$ such that if $s$ belongs to $\mathcal{B}_{2,\infty}^\delta(R)\cap W_\gamma(R')\cap\LL^{\infty}(R'')$ and satisfies
$$d^2(s,\calS_0)>
C(\alpha,\beta,R'',\delta,\gamma) \Bigg(\frac{\ln L}{L}+ R^2 \left(\frac{\ln L}{L}\right)^{2\delta}
+R'^2 \left(\frac{\ln L}{L}\right)^{\frac{2\gamma}{1+2\gamma}}\\+R'^{2+4\gamma}\left(\frac{\ln L}{L}\right)^{2\gamma}\Bigg),$$
then
$$ \PP_s\pa{\Phi_{\alpha}^{(2)}=0}\leq \beta. $$
In particular, when $\delta\geq \gamma/(1+2\gamma)$, there exist some positive constants $L_0(\delta,\gamma)$ and $C(\alpha,\beta,R,R',R'',\delta,\gamma)$ such that if $L>L_0(\delta,\gamma)$, then
$$\rho(\Phi_\alpha^{(2)},\mathcal{B}_{2,\infty}^\delta(R)\cap W_\gamma(R')\cap \LL^\infty(R''),\beta)\leq
C(\alpha,\beta,R,R',R'',\delta,\gamma)\left(\frac{\ln L}{L}\right)^{\frac{\gamma}{1+2\gamma}}.$$
\end{Prop}

\emph{Comments.}

\begin{enumerate}
\item Our second testing procedure is still adaptive: indeed, for large $L$, it achieves the lower bounds for the minimax separation
rates over all the spaces $\mathcal{B}_{2,\infty}^\delta(R)\cap W_\gamma(R')\cap \LL^\infty(R'')$ with ${\gamma}/(1+2\gamma)\leq \delta< \gamma/2$ simultaneously (see Theorem \ref{minorations}). In this case, we also remark that these rates are so large that there is no further price to pay for adaptivity in the sense that the upper bound does not involve any extra logarithmic factor. To our knowledge, this phenomenon is completely new for nonparametric testing procedures.
\item Our second procedure achieves the lower bounds for the minimax separation
rates over all the spaces $\mathcal{B}_{2,\infty}^\delta(R)\cap W_\gamma(R')\cap \LL^\infty(R'')$ with ${\gamma}/(1+2\gamma)\leq \delta <\gamma/2$ simultaneously, but it does not  when $\delta \geq \max{}(\gamma/2,\gamma/(1+2\gamma))$. To obtain a test that achieves the minimax separation
rates in both cases, our two procedures need to be combined.
\end{enumerate}

\subsection{Uniform separation rates of the combined procedure}

\begin{Cor}\label{vitesses3}

Assume that $\ln \ln L\geq 1$. Fix some level $\alpha$ and $\beta$ in $]0,1[$. Let $\Phi_{\alpha/2}^{(1)}$ be the level $\alpha/2$ test
 defined by (\ref{fonctiontest1}) with $\mathcal{J}=\{1,\ldots,\left\lfloor\log_2(L^2/(\ln\ln L)^3)\right \rfloor\}$ and $W_{J}=\ln|\mathcal{J}|$ for every $J$ in $\mathcal{J}$.  Let $\Phi_{\alpha/2}^{(2)}$ be the level $\alpha/2$ test defined by (\ref{fonctiontest2}) with $\bJ=\left\lfloor\log_2(L/\ln L)\right \rfloor$.
We consider $\Phi_\alpha^{(3)}=\sup{}
\{\Phi_{\alpha/2}^{(1)},\Phi_{\alpha/2}^{(2)}\}$.

$(i)$ For all $\delta>0$ and $\gamma>0$, there exist some positive constants $L_0(\delta)$ and $C(\alpha,\beta,R,R'',\delta)$ such that if $L>L_0(\delta)$, then
$$\rho(\Phi_\alpha^{(3)},\mathcal{B}_{2,\infty}^\delta(R)\cap W_\gamma(R')\cap\LL^{\infty}(R''),\beta)\leq
C(\alpha,\beta,R,R'',\delta) \left(\frac{\sqrt{\ln\ln L}}{L}\right)^{\frac{2\delta}{4\delta +1}}.$$

 $(ii)$ For all $(\delta,\gamma)$ such that $\gamma/(2\gamma+1)\leq \delta$, there exist some positive constants $L_0(\delta,\gamma)$ and $C(\alpha,\beta,R,R',R'',\delta,\gamma)$ such that if $L>L_0(\delta,\gamma)$, then
$$\rho(\Phi_\alpha^{(3)},\mathcal{B}_{2,\infty}^\delta(R)\cap W_\gamma(R')\cap\LL^{\infty}(R''),\beta)\leq
C(\alpha,\beta,R,R',R'',\delta,\gamma)\left(\frac{\ln L}{L}\right)^{\frac{\gamma}{1+2\gamma}}.$$

\end{Cor}

\emph{Comment.} Since $$\PP_s\pa{\Phi_{\alpha}^{(3)}=0}\leq \inf{}\left\{\PP_s\pa{\Phi_{\alpha/2}^{(1)}=0},\PP_s\pa{\Phi_{\alpha/2}^{(2)}=0}\right\},$$ the proof of this result directly comes from Proposition \ref{vitesses1} and Proposition \ref{vitesses2}.\\
This final procedure actually matches the lower bounds of Theorem \ref{minorations} and is consequently adaptive for the whole set of parameters $(\delta,\gamma)$ such that $\gamma/(2\gamma+1)\leq \delta$ (up to a  $\ln \ln L$ factor when $\delta \geq \gamma/2$). This also proves that the lower bounds of Theorem \ref{minorations} are sharp for this set of parameters.

\section{Simulation study}\label{sectionsimu}

We aim in this section at studying the performances of our tests from a practical point of view. We consider several intensities $s$ defined on $[0,1]$ such that
$\int_0^1 s(x) dx=1$. $N$ denotes here a Poisson process with intensity $L s$ on $[0,1]$ with respect to the Lebesgue measure, and $\PP_s$ the distribution of this process. We denote by $s_0$ the intensity which is constant (equal to 1) on $[0,1]$. We choose $L=100$ and a level of test $\alpha=0.05$.

Let us now recall that our first procedure may be based on the test statistics
 \beqe  \mathcal{T}_{\alpha}^{(1)}=\sup{J \in\mathcal{J}}\left(T_J'-q_{J}'^{(N_L)}(u_{\alpha}'^{(N_L)}/|\mathcal{J}|)\right), \eeqe
where $q_{J}'^{(n)}(u)$ denotes the $(1-u)$ quantile of $ T_J' | N_L=n $ under the hypothesis that $s=s_0$, and $u_{\alpha}'^{(n)}$ is chosen such that :
$$u_{\alpha}'^{(n)}=\sup{}\ac{u \in ]0,1[, \PP_{s_0}\pa{\sup{J\in\mathcal{J}}\left(
T_J'-q_J'^{(n)}(u/|\mathcal{J}|)\right)>0\Bigg|N_L=n} \leq \alpha }.$$
The null hypothesis $(H_0)$  "$s=s_0$" is rejected when $\mathcal{T}_{\alpha}^{(1)}>0$.

We choose $\mathcal{J}=\{1,\ldots,6\}$. For $40\leq n\leq 160$, we estimate the quantities $u_{\alpha}'^{(n)}$
and the quantiles $q_{J}'^{(n)}(u_{\alpha}'^{(n)}/|\mathcal{J}|)$ for all $J$ in $\mathcal{J}$.
These estimations are based on the simulation of $200 000$  independent samples with size $n$, uniformly distributed on [0,1]. Half of the samples is used to estimate the
 quantiles $ q_{J}'^{(n)}(u/|\mathcal{J}|)$ for $u$ varying on a grid over $[0,1]$, and the other samples are used to estimate the probabilities occurring in the definition of $u_{\alpha}'^{(n)}$. Finally, $u_{\alpha}'^{(n)}$ is estimated by the largest value on the grid such that  these estimated probabilities are smaller than $\alpha$.\\

Let us also recall that our second procedure is based on the test statistics
 \beqe  \mathcal{T}_{\alpha}^{(2)}=\sup{\lambda\in \Lambda_\bJ}
\left(T_\lambda-q_{\lambda}^{(N_L)}(u_{\lambda,\alpha}^{(N_L)})\right), \eeqe
where $q_{\lambda}^{(n)}(u)$ denotes the $(1-u)$ quantile of $ T_{\lambda} | N_L=n $ under the hypothesis that $s=s_0$. For  $\lambda=(j,k)\in \Lambda_{\bJ}$,
$ u_{\lambda,\alpha}^{(n)} = {u_{\alpha}^{(n)}}/\pa{ 2^j \bJ}$ with $u_{\alpha}^{(n)}$ defined by (\ref{ualpha}).
The null hypothesis $(H_0)$  "$s=s_0$" is rejected when $\mathcal{T}_{\alpha}^{(2)}>0$.\\

We choose $\bJ=6$. For $40\leq n\leq 160$, we estimate the quantities $u_{\alpha}^{(n)}$
and the quantiles $ q_{\lambda}^{(n)}( u_{\lambda,\alpha}^{(n)})$ for all $\lambda \in  \Lambda_{\bJ}$.
These estimations are based on the simulation of $200 000$  independent samples with size $n$, uniformly distributed on [0,1]. Half of the samples is used to estimate the
 quantiles $ q_{(j,k)}^{(n)} \pa{u/(2^j \bJ)}$ for $u$ varying on a grid over $[0,1]$, and the other samples are used to estimate the probabilities that occur in (\ref{ualpha}). Finally, we estimate $ u_{\alpha}^{(n)}$ in the same way as in the first procedure.\\

At this stage, we can estimate the powers of the two tests under various alternatives. The chosen alternatives are intensities that have already been studied  among others by Reynaud-Bouret and Rivoirard \cite{PatouVincent}, in the estimation problem. Since we are particularly interested in detecting the homogeneity of a Poisson process when the alternatives may be very irregular, we focus on the functions defined by:
\begin{eqnarray*}
&& s_1(x)=(1+\varepsilon) \1_{[0,0.125[}(x) + (1-\varepsilon) \1_{[0.125,0.25[}(x)+\1_{[0.25,1]}(x), \\
&& s_2(x)= \pa{1+\eta\sum_j \frac{h_j}{2}(1+\mbox{sgn}(x-p_j))} \frac{\1_{[0,1]}(x)}{C_2(\eta)},\\
&& s_3(x)= (1-\varepsilon) \1_{[0,1 ]}(x) + \varepsilon \pa{\sum_{j} g_j\pa{1+\frac{|x-p_j|}{w_j}}^{-4}}\frac{\1_{[0,1]}(x)}{0.284},\\
\end{eqnarray*}
where
$$
\left\{\begin{tabular}{ccccccccccccccc}
$p$&$=$&[& 0.1&0.13&0.15&0.23&0.25&0.4&0.44& 0.65&0.76& 0.78& 0.81& ]\\
$h$&$=$& [&4&-4&3&-3&5&-5&2&4&-4&2&-3&] \\
$g$&$=$&[&4&5&3&4&5&4.2&2.1&4.3&3.1&5.1&4.2&] \\
$w$&$=$&[&0.005&0.005&0.006&0.01&0.01&0.03&0.01&0.01&0.005&0.008&0.005&]
\end{tabular},\right.$$
$0<\varepsilon\leq 1$, $0<\eta\leq 2$, and $C_2(\eta)$ is such that $\int_0^1 s_2(x)dx=1$. \\

These alternatives, for particular values of the parameters, are represented in Figure 1.

In  Figure 2, we represent the histograms of one simulated sample for some of these alternatives  and for a constant intensity on $[0,1]$. Note that these histograms are clearly not sufficient to separate the alternatives from the null hypothesis.

\bigskip\bigskip

\begin{figure}[hbt]\label{alter}
\includegraphics[width=16cm,height=5cm]{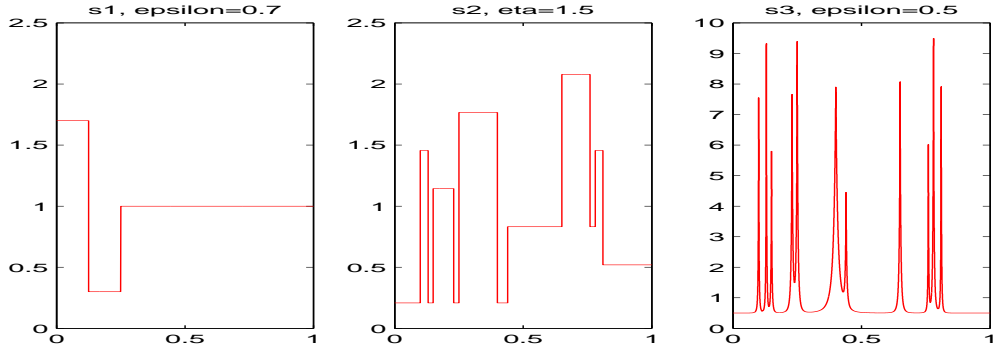}
\caption { Functions $s_1,s_2, s_3$}
\end{figure}

\newpage

\begin{figure}[hbt]\label{histo}
\includegraphics[width=16cm,height=9cm]{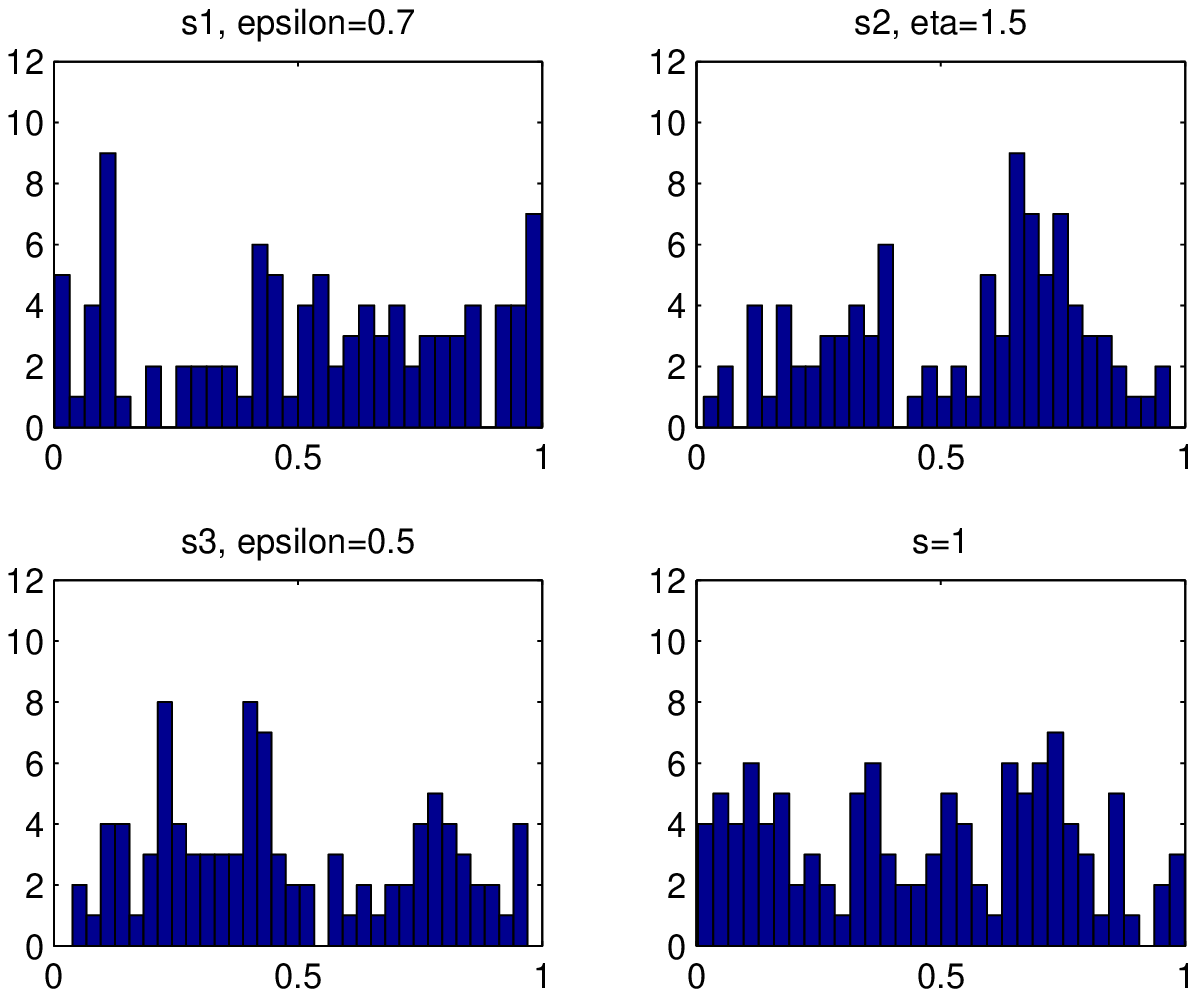}
\caption {Histograms of one simulated Poisson process}
\end{figure}

We also consider two monotonous alternatives defined by :
\begin{eqnarray*}
&& s_4(x)=(1-\varepsilon) \1_{[0,0.75[}(x) + (1+3\varepsilon) \1_{[0.75,1]}(x), \\
&& s_5(x)= (1-\varepsilon) \1_{[0,1 ]}(x) + \varepsilon \beta
x^{\beta-1}\1_{[0,1 ]}(x) ,
\end{eqnarray*}
where $0<\varepsilon<1$, and $\beta>1$.\\
These alternatives, for particular values of the parameters, are represented in Figure 3.

In  Figure 4, we represent the histograms of one simulated sample for some of these alternatives.

\bigskip

\begin{figure}[hbt]\label{altercroiss}
\centering
\includegraphics[width=16cm,height=5.5cm]{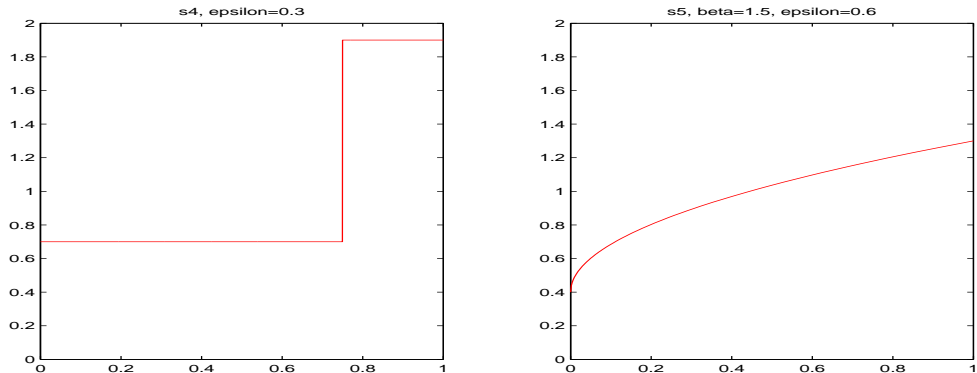}
\caption { Functions $s_4$ and $s_5$}
\end{figure}

\newpage

\begin{figure}[hbt]\label{histocroiss}
\centering
\includegraphics[width=16cm,height=4.5cm]{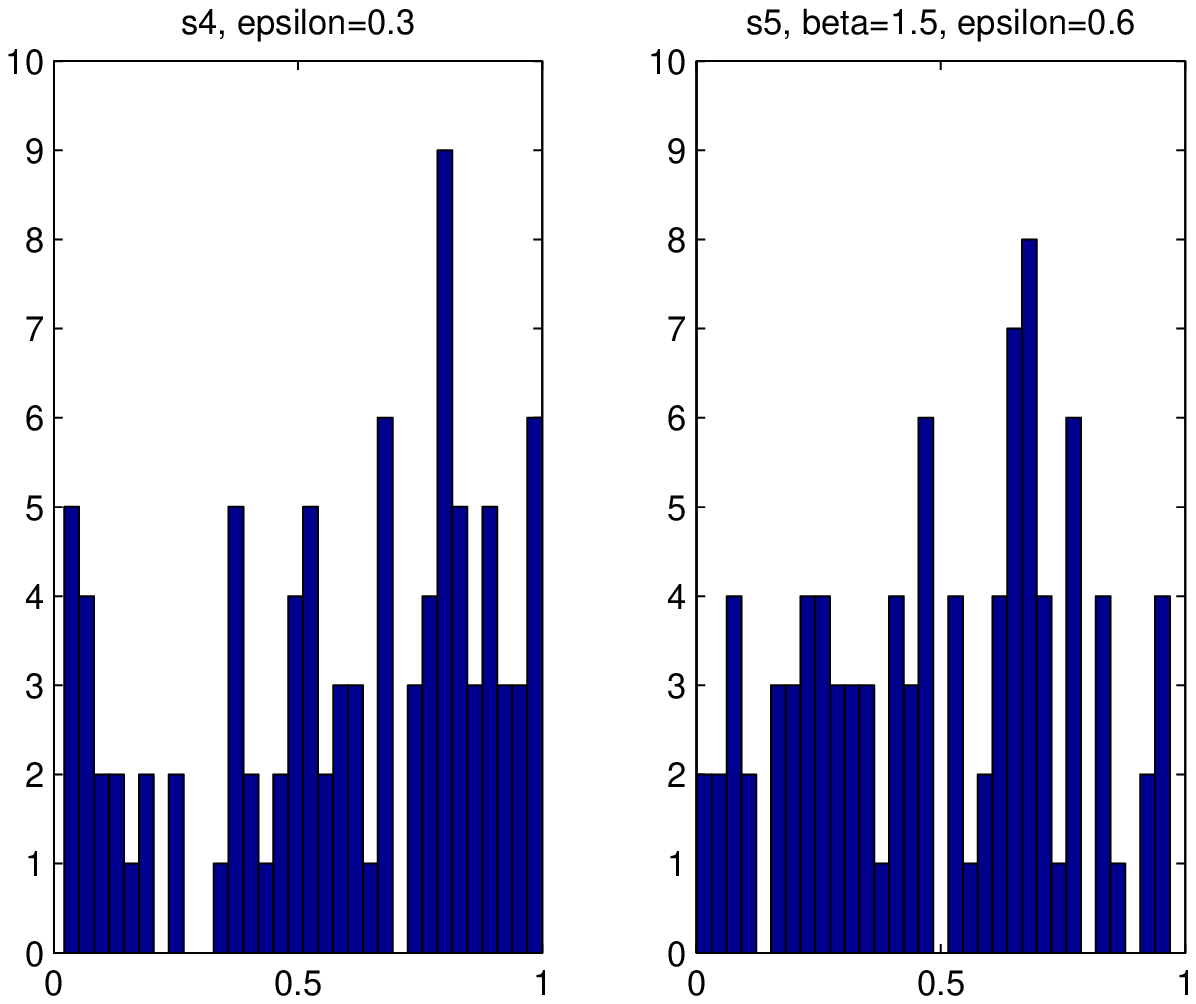}
\caption { Histograms of one simulated Poisson process}
\end{figure}

For each alternative $s$, we simulate $20000$ Poisson processes with intensity $Ls$ on $[0,1]$, and we estimate the powers of our two tests by :
$$ \hat{P}_1= \frac1{20000} \sum_{k=1}^{20000} \1_{\mathcal{T}_{\alpha}^{(1),k}>0},$$
and
$$\hat{P}_2= \frac1{20000} \sum_{k=1}^{20000} \1_{\mathcal{T}_{\alpha}^{(2),k}>0},$$
where $\mathcal{T}_{\alpha}^{(1),k }$ and $\mathcal{T}_{\alpha}^{(2),k }$ are the test statistics $\mathcal{T}_{\alpha}^{(1)}$ and $\mathcal{T}_{\alpha}^{(2)}$ computed for the $k$th simulated Poisson process.

We compare the obtained estimated powers with the estimated powers of the classical Kolmogorov and Smirnov's test applied to the Poisson process conditionally on the event "the number of points of the Poisson process is $n$". The estimated powers of Kolmogorov and Smirnov's test denoted by $\hat{P}_{KS}$ are also obtained by $20000$ simulations of a Poisson process with intensity $Ls$ on $[0,1]$.

The estimated powers are furthermore compared to the estimated powers of the tests studied in practice by the other authors. Such tests are in fact devoted to the particular case of increasing alternatives, which may be relevant in reliability contexts involving repairable systems. Bain, Engelhardt and Wright \cite{Bainetal} and Cohen and Sackrowitz \cite{Cohen} consider in these contexts six well known tests. They show that two of these six tests, namely the so-called Laplace and $Z$ tests (respectively studied first by Cox \cite{Cox} and Crow \cite{Crow}) are preferable to use.

The Laplace test is based on the statistics
$$\mathcal{T}_{\alpha}^{(La)}=\sum_{l=1}^{N_L} X_l-q_{La}^{(N_L)}(\alpha),$$
where $(X_1,\ldots,X_{N_L})$ are the points of the process, and for every $n$, $q_{La}^{(n)}(\alpha)$ is the $(1-\alpha)$ quantile of the sum of $n$ independent random variables uniformly distributed on $[0,1]$.

The $Z$ test is based on the statistics
$$\mathcal{T}_{\alpha}^{(Z)}=2\sum_{l=1}^{N_L} \ln(X_l)+q_{Z}^{(N_L)}(\alpha),$$
where for every $n$, $q_{Z}^{(n)}(\alpha)$ is the $\alpha$ quantile of the chi square distribution with $2n$ degrees of freedom.

Assuming that the intensity $s$ is increasing, the null hypothesis $(H_0)$  "$s$ is constant on $[0,1]$" is rejected when $\mathcal{T}_{\alpha}^{(La)}>0$ or $\mathcal{T}_{\alpha}^{(Z)}>0$.

The readers need to be aware that these tests are especially constructed to detect homogeneity against increasing trend, when reading the estimated power tables.\\

Let us now present the results we obtained for the different tests. The estimated powers for Poisson processes with intensities $L s_1$, $Ls_2$, $Ls_3$, $Ls_4$, and $Ls_5$ with various values of the parameters are given in the following tables.

\smallskip

\begin{center}
\underline{Alternatives $s_1$:}
\end{center}

\begin{center}
\begin{tabular}{|c|c|c|c|c|c|c|c|c|}
\hline
$\varepsilon$  & $0$ &$0.5$& $0.6$ &$0.7$ & $0.8$ & $0.9$  & $1$\\
\hline
\hline
{\bf $\hat{P}_1$} &$0.05$& $0.25$& $0.39$& $0.56$ & $0.73$ & $0.89$ & $0.98$\\
\hline
{\bf $\hat{P}_2$} &$0.05$& $0.33 $& $0.52 $& $0.72$ & $0.87$ & $0.96$ & $1$\\
\hline
{\bf $\hat{P}_{KS}$} &$0.05$& $0.09 $& $0.13 $& $0.19$ & $0.27$ & $0.37$ & $0.48$ \\
\hline
{\bf $\hat{P}_{La}$} &$0.05$& $0.03 $& $0.03 $& $0.03$ & $0.03$ & $0.02$ & $0.02$ \\
\hline
{\bf $\hat{P}_{Z}$} &$0.05$& $0.01 $& $0.01$& $0.01$ & $0.01$ & $0.01$ & $0.01$ \\
\hline
\end{tabular}
\end{center}

\begin{center}
\underline{Alternatives $s_2$:}
\end{center}

\begin{center}
\begin{tabular}{|c|c|c|c|c|c|}
\hline
$\eta$ &$0$& $0.5$&   $ 1 $& $1.5$ &$2$ \\
\hline
$C_2(\eta)$&$ 1 $ & $ 2.27$ & $ 3.54 $ & $4.81 $ & $6.08$ \\
\hline
\hline
{\bf $\hat{P}_1$} &$0.05  $ &$0.61$ & $0.87$   & $0.94$& $0.97$ \\
\hline
{\bf $\hat{P}_2$} &$0.05  $& $ 0.41$ & $0.64$   & $0.75$ & $0.80$\\
\hline
{\bf $\hat{P}_{KS}$} &$0.05  $&$ 0.14$ & $0.25$   & $0.34$& $0.39$\\
\hline
{\bf $\hat{P}_{La}$} &$0.05  $&$ 0.05$ & $0.06 $   & $0.06$& $0.06$\\
\hline
{\bf $\hat{P}_{Z}$} &$0.05  $&$ 0.26$ & $0.39$   & $0.46$ & $0.51$\\
\hline
\end{tabular}
\end{center}

\begin{center}
\underline{Alternatives $s_3$:}
\end{center}

\begin{center}
\begin{tabular}{|c|c|c|c|c|c|c|}
\hline
$\varepsilon$  & $0$ & $0.2$ & $0.3$& $0.4$ &$0.5$ & $0.6$ \\
\hline
\hline
{\bf $\hat{P}_1$} &$0.05$& $0.28$ & $ 0.65 $& $0.91$ & $0.99 $& $1$ \\
\hline
{\bf $\hat{P}_2$} &$0.05$& $0.20$ & $ 0.43 $& $0.71$ & $0.90 $& $0.98$ \\
\hline
{\bf $\hat{P}_{KS}$} &$0.05$& $0.11$&  $0.21$& $0.37$& $0.56$& $0.76$ \\
\hline
{\bf $\hat{P}_{La}$} &$0.05$& $0.01$&  $0.00$& $0.00$& $0.00$& $0.00$ \\
\hline
{\bf $\hat{P}_{Z}$} &$0.05$& $0.02$&  $0.02$& $0.02$& $0.02$& $0.01$ \\
\hline
\end{tabular}
\end{center}

\begin{center}
\underline{Alternatives $s_4$:}
\end{center}

\begin{center}
\begin{tabular}{|c|c|c|c|c|c|}
\hline
$\varepsilon$  & $0$ & $0.1$ & $0.2$ & $0.3$& $0.4$\\
\hline
\hline
{\bf $\hat{P}_1$} &$0.05$& $ 0.20$ & $0.69$ & $0.97$ & $1$\\
\hline
{\bf $\hat{P}_2$} &$0.05$ & $0.17$ & $0.62$   &$0.95$ & $1$\\
\hline
{\bf $\hat{P}_{KS}$} &$0.05$ & $0.26$ &$0.77$&   $0.98$ & $1$\\
\hline
{\bf $\hat{P}_{La}$} &$0.05$ &$0.37$ & $0.82$ &$0.98$ & $1$\\
\hline
{\bf $\hat{P}_{Z}$} &$0.05$ &   $0.24$ &$   0.57$ &$    0.85$ &$   0.97$\\
\hline
\end{tabular}
\end{center}

\begin{center}
\underline{Alternatives $s_5$:}
\end{center}

\begin{center}
\begin{tabular}{|c|c|c|c|c|c|c|}
\hline
$(\beta,\varepsilon)$  & $(1.5,0.2)$ & $(1.5,0.6)$ & $(1.5,1)$ & $(2,0.2)$ & $(2,0.6)$& $(2,1)$\\
\hline
\hline
{\bf $\hat{P}_1$} &$0.20$& $0.49$ & $0.79$& $0.24$ & $0.62$& $1$ \\
\hline
{\bf $\hat{P}_2$} &$0.18$& $0.43$ & $0.69$& $0.24$ & $0.62$& $1$ \\
\hline
{\bf $\hat{P}_{KS}$} &$0.22$& $0.56$&  $0.91$& $0.24$& $0.62$& $1$ \\
\hline
{\bf $\hat{P}_{La}$} &$0.24$& $0.60$&  $0.98$& $0.24$& $0.62$& $1$ \\
\hline
{\bf $\hat{P}_{Z}$} &$0.24$&  $0.61$& $0.99$& $0.24$& $0.62$ & $1$ \\
\hline
\end{tabular}
\end{center}

\emph{Comments.}

\begin{enumerate}
\item It first emerges from these results that when the alternatives are not increasing, our two tests have estimated powers significantly larger than the Laplace and $Z$ tests that are designed for increasing alternatives, but also than Kolmogorov and Smirnov's test. Furthermore, we can not give prior arguments to choose one of our two tests rather than the other one in these case. Indeed, we can notice that the first one is more powerful than the second one for alternatives $s_2$ which are rather smooth, but also for alternatives $s_3$ which are very irregular.
    Thus, in the case of non increasing alternatives such as $s_1$, $s_2$ and $s_3$, or in the practical situations of our interest such as the study of occurrences on DNA sequences where the intensities may have some localized spikes, this should argue in favor of the choice of our combined procedure.

\item As for the increasing alternatives, the specific Laplace and $Z$ tests remain as expected the most powerful ones, except for the alternatives $s_4$, that are not as smooth as the $s_5$ alternatives. Kolmogorov and Smirnov's test is also often more powerful than our tests. However, we know that in the case of smooth alternatives, we could probably significantly improve the estimated powers of our first test by using the Fourier basis instead of the Haar basis. Since our first test is very similar to Fromont and Laurent's \cite{Beamag} one in the density model, we  refer to this paper for more details. We could also consider a new test combining for instance our first test with the Laplace test.
\end{enumerate}

\section{Proofs}

\subsection{Proof of Theorem \ref{minorations}}

Since it is easier to argue in terms of errors of second kind than in terms of minimax separation rates directly, we start by defining for all $\calS\subset\LL^2([0,1])$,
$$\beta(\calS)=\inf{\Phi_\alpha}\sup{s\in\calS} \PP_s(\Phi_\alpha=0),$$
where the infimum is taken over all level $\alpha$ tests $\Phi_\alpha$ and by stating a useful and well-known lemma.
\begin{Lemma}\label{errdeux}
Let $r$ be a positive number, and $\calS$, $\calS'$ be subsets of $\LL^2([0,1])$.\\
 $(i)$ If $\beta(\{s\in\calS,d(s,\calS_0)\geq r\})\geq \beta,$
 then
$$\underline{\rho}(\calS, \alpha,\beta)\geq r.$$
$(ii)$ If $\calS'\subset \calS$, then $\beta(\calS)\geq \beta(\calS')$.
\end{Lemma}
The proof of the lemma is straightforward.

Our aim here is to construct finite sets $S_{M,D,r}$ such that
\begin{equation}\label{condition1}
S_{M,D,r}\subset \{s\in \mathcal{B}_{2,\infty}^{\delta}(R)\cap W_{\gamma}(R')\cap \mathbb{L}^{\infty}(R''),d(s,\calS_0)\geq r\},
\end{equation}
 and that
\begin{equation}\label{condition2}
\beta(S_{M,D,r})\geq \beta,
\end{equation} with $r$ as large as possible.

These finite sets are based on a family of functions $\{\p_{M,i}, i \in \{1,...,M\}\}$ such that for all $x \in [0,1]$, $\p_{M,i}(x)=\p(Mx-i+1),$ where $\p$ is a function on $[0,1]$  such that
\begin{eqnarray}\label{condvarphi}
\int_0^1\p(x) dx =0 \nonumber\\
\int_0^1 \p(x)^2dx =1 \\
\forall x \in [0,1], |\p(x)|\leq \rho.\nonumber
\end{eqnarray}

For $r>0 $, and $D\leq M$, we introduce the set
\begin{equation}\label{defreseau}
S_{M,D,r}=\ac{ s_{\xi,\Delta,r}= \rho \1_{[0,1]} +r \sqrt{\frac MD}\sum_{i=1}^M \Delta_i \xi_i \p_{M,i},\  \xi \in \{-1,+1\}^M, \Delta \in \{0,1\}^M, \sum_{i=1}^M \Delta_i=D}.
\end{equation}

As a first step, we notice that the functions $s_{ \xi,\Delta,r}$'s are positive as soon as $ r^2\leq D/M$ and that for every $s_{ \xi,\Delta,r} \in S_{M,D,r} $, $d(s_{\xi,\Delta,r},\calS_0)^2=\norm{s_{ \xi,\Delta,r}-\rho\1_{[0,1]}}^2=r^2$ (see (\ref{condvarphi})).

As a second step, we want to find which positive $r$ leads to $\beta(S_{M,D,r})\geq \beta$.

Let us recall a fundamental lemma which can be found in Ingster~\cite{Ing} or Baraud~\cite{yannick} for  other frameworks.
\begin{Lemma}\label{ingyannick}
Let $\nu$ be a probability measure on $S_{M,D,r}$ and let $\sigma\sim \nu$. Let $\PP_{\nu}$ be the distribution of a point process $N$ such that the conditional distribution of $N$ given that $\sigma=s$ is a Poisson process with intensity $s$. Let $\PP_0$ be the distribution of a Poisson process with constant intensity given by $\rho\1_{[0,1]}$, and $\EE_{0}$ denote the expectation with respect to $\PP_0$.
Let $L_{\nu}$ be the likelihood ratio $L_{\nu}= d\PP_{\nu}/d\PP_0.$ Then
$$\beta(S_{M,D,r})\geq 1-\alpha-\frac{1}{2}\left(\EE_0[L_\nu^2(N)]-1\right)^{1/2}.$$
\end{Lemma}
\begin{proof}
The proof is obtained by rather straightforward computations. One has
\begin{eqnarray*}
\beta(S_{M,D,r}) &\geq& \inf{\Phi_\alpha} \int \PP_s(\Phi_\alpha=0)d\nu(s) \\
&\geq & 1 - \sup{\Phi_\alpha} \int \PP_s(\Phi_\alpha=1)d\nu(s) \\
&\geq & 1-\sup{\Phi_\alpha} \left[|\PP_\nu(\Phi_\alpha=1)-\PP_{0}(\Phi_\alpha=1)|+|\PP_{0}(\Phi_\alpha=1)|\right]\\
&\geq & 1-\alpha-\norm{\PP_\nu-\PP_{0}}_{TV},
\end{eqnarray*}
where $\norm{.}_{TV}$ corresponds to the total variation norm.
Hence,
\begin{eqnarray*}
\beta(S_{M,D,r}) &\geq& 1 -\alpha- \frac12 \EE_{0}\left[\left|\frac{d\PP_\nu}{d\PP_{0}}-1\right|\right]\\
&\geq &  1 -\alpha- \frac12\EE_{0}\left[\left|L_\nu(N)-1\right|^2\right]^{1/2}.
\end{eqnarray*}
But $\EE_{0}[L_\nu(N)]=1$. So $\beta(S_{M,D,r}) \geq 1 -\alpha-\left(\EE_0[L_\nu^2(N)]-1\right)^{1/2}/2$.
\end{proof}

Regarding Lemma \ref{ingyannick}, we still have to find a distribution $\nu$ and $r$ such that $\EE_{0}[L_{\nu}^2(N)]\leq 1 +4(1-\alpha-\beta)^2$ which implies that $\beta(S_{M,D,r})\geq\beta.$

Let $\xi=(\xi_1,\ldots,\xi_M)$ be a random vector, such that the $\xi_i$'s are i.i.d. Rademacher variables, taking the values $+1$ and $-1$ with probability $1/2$. Let $\Delta=(\Delta_1,\ldots,\Delta_M)$ be a random vector, independent of $\xi$ and defined by $\Delta_i=\1_{i\in\mathcal{I}}$, where $\mathcal{I}$ is a set of $D$ indices drawn at random from $\{1,\ldots,M\}$ without replacement.

Then the random function $s_{\xi,\Delta,r}= \rho \1_{[0,1]} +r \sqrt{\frac MD}\sum_{i=1}^M \Delta_i \xi_i \p_{M,i}$ belongs to $S_{M,D,r}$, which allows to take its distribution as $\nu$.

Let us denote by $\EE_\xi$ the expectation with respect to the variable $\xi$ and by $\EE_\mathcal{I}$ the expectation with respect to the
random set $\mathcal{I}$ defined above. By definition,
$L_\nu=\int d\PP_s/d\PP_{0} d\nu(s).$ Hence
$L_\nu(N)=\EE_\mathcal{I}\EE_{\xi}\left[\exp\left(\int_0^1
\ln(s_{ \xi,\Delta,r}(x)/\rho) dN_x\right)\right].$
This can be rewritten as
$$
L_\nu(N)=\EE_\mathcal{I}\EE_{ \xi}\left[\prod_{i=1}^M
\exp \left(\int_{\left]\frac{i-1}{M},\frac iM\right]} \ln\left(1+r\sqrt{\frac MD}  \xi_i \Delta_i
\frac{\p_{M,i}(x)}{\rho}\right) dN_x\right)\right]
=\EE_\mathcal{I}\left[\prod_{i\in \mathcal{I}} A_i\right]
$$
where
\begin{multline*}
A_i = \frac12 \Bigg(\exp\left(\int_{\left]\frac{i-1}{M},\frac iM\right]} \ln\left(1+r\sqrt{\frac MD}
\frac{\p_{M,i}(x)}{\rho}\right) dN_x\right)+\\
\exp\left(\int_{\left]\frac{i-1}{M},\frac iM\right]} \ln\left(1-r\sqrt{\frac MD}
\frac{\p_{M,i}(x)}{\rho}\right) dN_x\right)\Bigg).
\end{multline*}

Let $\mathcal{I}'$ be a random set of indices with the same distribution as $\mathcal{I}$ and independent of $\mathcal{I}$. Then,

\begin{eqnarray*}
\EE_{0}[L_\nu^2(N)] &=& \EE_{0}\EE_\mathcal{I}\EE_{\mathcal{I}'}\left[\prod_{i\in \mathcal{I}}A_i\prod_{i\in \mathcal{I}'}A_i\right]\\
&=&\EE_{0}\EE_\mathcal{I}\EE_{\mathcal{I}'}[\prod_{i\in \mathcal{I}\setminus\mathcal{I}' }A_i\prod_{i\in \mathcal{I}'\setminus\mathcal{I}}A_i\prod_{i\in \mathcal{I}\cap\mathcal{I}'}A_i^2 ]
\end{eqnarray*}
But under the distribution $\PP_0$, the variables $A_i$'s are mutually independent since they only depend on the integrals of the Poisson process on intervals with disjoint support. Consequently,
\begin{equation}\label{esperanceLcarre}
\EE_{0}[L_\nu^2(N)]=\EE_\mathcal{I}\EE_{\mathcal{I}'}\left[\prod_{i\in \mathcal{I}\setminus\mathcal{I}' }\EE_{0}[A_i]\prod_{i\in \mathcal{I}'\setminus\mathcal{I}}\EE_{0}[A_i]\prod_{i\in \mathcal{I}\cap\mathcal{I}'}\EE_{0}[A_i^2]\right].
\end{equation}

We now need to compute $\EE_{0}[A_i]$ and $\EE_{0}[A_i^2]$, and we use the following lemma.

\begin{Lemma}\label{calculsprocessus}
Let $f$ be a function on $[0,1]$. Then with the above notations,
$$\EE_0\left[\exp\left(\int_{\left]\frac{i-1}{M},\frac iM\right]} f(x) dN_x\right)\right]=\exp\left(\int_{\left]\frac{i-1}{M},\frac iM\right]}\left(\exp(f(x))-1\right)\rho L dx\right).$$
\end{Lemma}

\begin{proof}
When $N$ has the constant intensity $\rho \1_{[0,1]}$, we know that conditionally on the event "the number of points $N_{M,i}=N\left(\left]\frac{i-1}{M},\frac iM\right]\right)$ falling into $\left]\frac{i-1}{M},\frac iM\right]$ is $n$", the points of the process $X_1,\ldots,X_{N_{M,i}}$ in $\left]\frac{i-1}{M},\frac iM\right]$ obey the same law as a $n$-sample with uniform distribution on $\left]\frac{i-1}{M},\frac iM\right]$. Then, one can easily see that
\begin{eqnarray*}
\EE_0\left[\exp\left(\int_{\left]\frac{i-1}{M},\frac iM\right]} f(x) dN_x\right)\right]&=& \EE_0\left[\exp\left(\sum_{l=1}^{N_{M,i}} f(X_l)\right)\right]\\
&=&\EE_0\left[ \prod_{l=1}^{N_{M,i}}\EE\left[\exp(f(X_l))\big|N_{M,i} \right]\right]\\
&=&\EE_0\left[\exp\left(N_{M,i} \ln\left( \int_{\left]\frac{i-1}{M},\frac iM\right]} \exp(f(x)) M dx\right)\right)\right].
\end{eqnarray*}
Under $\PP_0$, $N_{M,i}$ has a Poisson distribution with parameter $\rho L/M$, therefore,
$$\EE_0\left[\exp\left(N_{M,i} \ln\left( \int_{\left]\frac{i-1}{M},\frac iM\right]} \exp(f(x)) M dx\right)\right)\right]=\exp\left(\frac{\rho L}{M}\left(\int_{\left]\frac{i-1}{M},\frac iM\right]}\exp(f(x))Mdx-1 \right)\right).$$
This concludes the proof.
\end{proof}

From Lemma \ref{calculsprocessus} and (\ref{condvarphi}), one has that
\begin{eqnarray*}
\EE_{0}[A_i]&=&\frac{1}{2}\exp\left(\int_{\left]\frac{i-1}{M},\frac iM\right]}\left(r\sqrt{\frac MD}
\frac{\p_{M,i}(x)}{\rho}\right)\rho L dx\right)\\
&&+\frac{1}{2}\exp\left(\int_{\left]\frac{i-1}{M},\frac iM\right]}\left(-r\sqrt{\frac MD}
\frac{\p_{M,i}(x)}{\rho}\right)\rho L dx\right)\\
&=&1.
\end{eqnarray*}
Moreover,
\begin{eqnarray*}
\EE_0[A_i^2]&=&\frac{1}{4}\EE_0\left[\exp\left(\int_{\left]\frac{i-1}{M},\frac iM\right]} 2\ln\left(1+r\sqrt{\frac MD}
\frac{\p_{M,i}(x)}{\rho}\right) dN_x\right)\right]\\
&&+\frac{1}{4}\EE_0\left[\exp\left(\int_{\left]\frac{i-1}{M},\frac iM\right]} 2\ln\left(1-r\sqrt{\frac MD}
\frac{\p_{M,i}(x)}{\rho}\right) dN_x\right)\right]\\
&&+\frac 12 \EE_0\left[\exp\left(\int_{\left]\frac{i-1}{M},\frac iM\right]} \ln\left(1-r^2\frac MD
\frac{\p_{M,i}^2(x)}{\rho^2}\right) dN_x\right)   \right].
\end{eqnarray*}
Using Lemma \ref{calculsprocessus} and (\ref{condvarphi}) again, we finally obtain that
$$\EE_0[A_i^2]=\cosh\left(\frac{r^2L}{\rho D}\right).$$

Hence, equation (\ref{esperanceLcarre}) gives

\begin{eqnarray*}
\EE_{0}[L_\nu^2(N)] &=& \EE_\mathcal{I}\EE_{\mathcal{I}'}\left[\prod_{i\in \mathcal{I}\cap\mathcal{I}'}\cosh\left(\frac{r^2L}{\rho D}\right)\right]\\
&=& \EE_\mathcal{I}\EE_{\mathcal{I}'}\left[\exp\left(| \mathcal{I}\cap\mathcal{I}'|\ln\cosh\left(\frac{r^2L}{\rho D}\right)\right)\right].
\end{eqnarray*}
For fixed $\mathcal{I}$, $| \mathcal{I}\cap\mathcal{I}'|$ is an hypergeometric variable with parameters $(M,D, D/M)$. Hence, we know from Aldous \cite{Aldous} p. 173,  that there exists a binomial variable $B$ with parameter $(D,D/M)$ such that $\EE_{\mathcal{I}'}\big[B\big|| \mathcal{I}\cap\mathcal{I}'|\big]=| \mathcal{I}\cap\mathcal{I}'|$. By Jensen's inequality, we obtain that
$$\EE_{0}[L_\nu^2(N)]\leq \EE_\mathcal{I}\EE_{\mathcal{I}'}\left[\exp\left(B\ln\cosh\left(\frac{r^2L}{\rho D}\right)\right)\right].$$
Setting $B=\sum_{i=1}^D B_i$ where the $B_i$'s are independent random Bernoulli variables with parameter $D/M$, we easily obtain that
\begin{equation}\label{majespLcarre}
\EE_{0}[L_\nu^2(N)]\leq \exp\left(D\ln\left(1+\frac{D}{M}\left(
\cosh\left(\frac{r^2L}{\rho D}\right)-1\right)\right)\right).
\end{equation}

From equation (\ref{majespLcarre}) and Lemma \ref{ingyannick}, we see that if
$$\exp \left(D \ln\left(1+\frac{D}{M} \left(\cosh\left(\frac{r^2L}{\rho D}\right)-1\right)\right)\right)\leq 1+4(1-\alpha-\beta)^2,$$
then $\beta(S_{M,D,r}) \geq \beta$.

Following Baraud's idea \cite{yannick} and setting $c=1+4(1-\alpha-\beta)^2$, since the function $\cosh$ is increasing on $[0,+\infty[$, we have that if
$$r^2\leq  \frac{\rho D}{L}\ln\left(1+\frac{M}{D^2}\ln c+\sqrt{2\frac{M}{D^2}\ln c+\left(\frac{M}{D^2}\ln c\right)^2}\right),$$
then
\begin{eqnarray*}
\cosh\left(\frac{r^2L}{\rho D}\right)-1&\leq &\frac{1}{2}\left(\frac{M}{D^2}\ln c+\sqrt{2\frac{M}{D^2}\ln c+\left(\frac{M}{D^2}\ln c\right)^2}-1\right)\\
&&+\frac{1}{2}\pa{\frac{M}{D^2}\ln c+\sqrt{2\frac{M}{D^2}\ln c+\left(\frac{M}{D^2}\ln c\right)^2}+1}^{-1}\\
&\leq &\frac{M\ln c}{D^2}.
\end{eqnarray*}
Hence
\begin{eqnarray*}
\exp \left(D \ln\left(1+\frac{D}{M} \left(\cosh\left(\frac{r^2L}{\rho D}\right)-1\right)\right)\right)&\leq &\exp\left(D\ln\left(1+\frac{\ln c}{D} \right)\right)\\
&\leq & c,
\end{eqnarray*}
and $\beta(S_{M,D,r}) \geq \beta$.
As a conclusion, we obtain the following result, where the second part of the Proposition comes from a direct computation (see Baraud \cite{yannick} for further details).
\begin{Prop}\label{erreur2surreseau}
Let $c=1+4(1-\alpha-\beta)^2$ and $S_{M,D,r}$ be the finite set defined by (\ref{defreseau}).
If
\begin{equation}\label{conditionr2}
r^2\leq \frac DM~\mbox{ and }r^2 \leq \frac{\rho D}{L}\ln\left(1+\frac{M}{D^2}\ln c+\sqrt{2\frac{M}{D^2}\ln c+\left(\frac{M}{D^2}\ln c\right)^2}\right),
\end{equation}
then $S_{M,D,r}\subset\{s, s\geq 0, d(s,\calS_0)=r\}$ and $\beta(S_{M,D,r})\geq \beta$.

If $\alpha+\beta\leq 0.59$ and
$$r^2\leq \frac DM \wedge \left[\frac{\rho D}{L}\ln\left(1+\frac{M}{D^2}\vee \sqrt{\frac{M}{D^2}}\right)\right],$$
then $S_{M,D,r}\subset\{s, s\geq 0, d(s,\calS_0)=r\}$ and $\beta(S_{M,D,r})\geq \beta$.
\end{Prop}

As a third step, we are now in position to find some $r$ (as large as possible) such that $S_{M,D,r}\subset \{s\in \mathcal{B}_{2,\infty}^{\delta}(R)\cap W_{\gamma}(R')\cap \mathbb{L}^{\infty}(R''),d(s,\calS_0)\geq r\}$ and that
$\beta(S_{M,D,r})\geq \beta$.\\

Let us consider the set $S_{M,D,r}$ defined by (\ref{defreseau}) with $\p=\1_{[0,1/2[}-\1_{[1/2,1[}$, $M=2^{J}$ and $\rho=1$.

Let $s \in S_{M,D,r}$, then $s$ can be rewritten as $s=\alpha_0\phi_0+\sum_{j \in \mathbb{N}}\sum_{k=0}^{2^j-1}\alpha_{(j,k)}\phi_{(j,k)}$, with
$$\alpha_0=1, \alpha_{(j,k)}=0 \mbox{ if }j \neq J,  \alpha_{(J,k)}^2=\frac{r^2}{D}\Delta_{k+1} \mbox{ for } \ k=0, \ldots 2^{J}-1.$$
Since $\sum_{k=1}^{2^{J}}\Delta_k=D$, the condition $r^2 \leq R^2 2^{-2J\delta}$ ensures that $S_{M,D,r} \subset  \mathcal{B}_{2,\infty}^{\delta}(R).$\\
Let us define, for all $t>0$,
$$ H(t)=\sum_{k=1}^{2^{J}} \frac{r^2}{D}\Delta_k \1_{\frac{r^2}{D}\Delta_k \leq t}.$$
In order to ensure that  $s$ belongs to $ W_{\gamma}(R')$, the function $H$ has to satisfy
$$\forall t>0, H(t) \leq {R'}^2 t^{\frac{2\gamma}{1+2\gamma}}.$$
Note that
$$H(t)=0 \mbox{ for } t<\frac{r^2}{D} \mbox{ and  } H(t)=H\left(\frac{r^2}{D}\right) \mbox{ for } t \geq \frac{r^2}{D}.$$
Hence, we only need to have that
$$H\left(\frac{r^2}{D}\right)\leq {R'}^2 \pa{\frac{r^2}{D}}^{\frac{2\gamma}{1+2\gamma}},$$
 which is equivalent to
$$r^2 \leq {R'}^{2(1+2\gamma)}D^{-2\gamma}.$$
Moreover, if $r^2\leq D/M$, then $\norm{s}_\infty\leq 2$. Hence when $R''\geq 2$, the condition
\begin{equation}\label{inclusion}
r^2 \leq \frac{D}{M}\wedge R^2 M^{-2\delta}\wedge {R'}^{2(1+2\gamma)}D^{-2\gamma}
\end{equation}
ensures that $S_{M,D,r} \subset
\mathcal{B}_{2,\infty}^{\delta}(R) \cap W_{\gamma}(R')\cap \mathbb{L}^{\infty}(R'')$.
From Proposition \ref{erreur2surreseau}, we can conclude that when $R''\geq 2$ and $\alpha+\beta\leq 0.59$,
if
\begin{equation}\label{condition3}
r^2 \leq \frac{D}{M}\wedge  \left(\frac{\rho D}{L}\ln\left(1+\frac{M}{D^2}\vee \sqrt{\frac{M}{D^2}}\right)\right)\wedge R^2 M^{-2\delta}\wedge {R'}^{2(1+2\gamma)}D^{-2\gamma},
\end{equation}
then (\ref{condition1}) and (\ref{condition2}) are both satisfied.

We now consider several cases, that are represented on the following figure.\\

\begin{center}

\scalebox{0.75}{\begin{picture}(0,0)%
\includegraphics{cas.pstex}%
\end{picture}%
\setlength{\unitlength}{3947sp}%
\begingroup\makeatletter\ifx\SetFigFont\undefined%
\gdef\SetFigFont#1#2#3#4#5{%
  \reset@font\fontsize{#1}{#2pt}%
  \fontfamily{#3}\fontseries{#4}\fontshape{#5}%
  \selectfont}%
\fi\endgroup%
\begin{picture}(8022,6190)(886,-11830)
\put(8626,-6136){\makebox(0,0)[lb]{\smash{{\SetFigFont{12}{14.4}{\rmdefault}{\mddefault}{\updefault}{$\delta=\gamma$}%
}}}}
\put(8701,-8836){\makebox(0,0)[lb]{\smash{{\SetFigFont{12}{14.4}{\rmdefault}{\mddefault}{\updefault}{$\delta=\gamma/2$}%
}}}}
\put(8551,-10486){\makebox(0,0)[lb]{\smash{{\SetFigFont{12}{14.4}{\rmdefault}{\mddefault}{\updefault}{$\delta=\frac{\gamma}{1+2\gamma}$}%
}}}}
\put(8776,-11761){\makebox(0,0)[lb]{\smash{{\SetFigFont{12}{14.4}{\rmdefault}{\mddefault}{\updefault}{$\gamma$}%
}}}}
\put(901,-5836){\makebox(0,0)[lb]{\smash{{\SetFigFont{12}{14.4}{\rmdefault}{\mddefault}{\updefault}{$\delta$}%
}}}}
\end{picture}%
}

\end{center}

In the following of this proof, $C$ will denote a positive constant that may depend on $\alpha,\beta,R,R',R'',\delta,\gamma$, and that may vary from one line to another.\\

{\bf Case 1.} If  $\delta<\gamma/2$, and $\delta \geq  \gamma/(1+2\gamma)$, we set
$$D=\left\lfloor\pa{\frac{L}{\ln L}}^{\frac{1}{1+2\gamma}}\right\rfloor,$$
 and
$$M=2^{J},\mbox { with } J= \left\lfloor\log_2\pa{\frac{L}{\ln L}}^{\frac{\gamma/\delta}{1+2\gamma}}\right \rfloor+1.$$

We first check that $D\leq M$ for $L$ large enough since $\delta\leq \gamma$.

Then,
$$R'^{2(1+2\gamma)} D^{-2\gamma}\geq C(L/\ln L)^{\frac{-2\gamma}{1+2\gamma}},$$
and
$$R^2M^{-2\delta}\geq C
(L/\ln L)^{\frac{-2\gamma}{1+2\gamma}}.$$
Finally, since
$$\frac{M}{D^2}\geq  (L/\ln L)^{\frac{\gamma/\delta-2}{1+2\gamma}}\to_{L\to+\infty} +\infty \mbox{ when } \gamma>2\delta,$$
then
$$\frac DL \ln\left(1+\frac M {D^2}\vee \sqrt{\frac M {D^2}}\right)\geq C
(L/\ln L)^{\frac{-2\gamma}{1+2\gamma}},$$
and
$$\frac DM \geq  C(L/\ln L)^{\frac{1-\gamma/\delta}{1+2\gamma}}\geq  C
(L/\ln L)^{\frac{-2\gamma}{1+2\gamma}}
\mbox{ for }L \mbox{ large enough.}$$

{\bf Case 2.} If $\gamma >1/2$ and $\delta \leq  \gamma/(1+2\gamma)$, one chooses
$$D=\left\lfloor\pa{\frac{L}{\ln L}}^{\frac{1}{1+2\gamma}}\right \rfloor,$$
and $$M= 2^{J} \mbox{ with } J=\left \lfloor\log_2( L/\ln L)\right \rfloor+1.$$

We first check that $D\leq M$ for $L$ large enough since $\gamma>0$.
Then,
$$R'^{2(1+2\gamma)}D^{-2\gamma}\geq C(L/\ln L)^{\frac{-2\gamma}{1+2\gamma}},$$
and
$$R^2M^{-2\delta}\geq C(L/\ln L)^{-2\delta}\geq
C(L/\ln L)^{\frac{-2\gamma}{1+2\gamma}}.$$
Since moreover
$$\frac{M}{D^2}\geq  (L/\ln L)^{1-\frac{2}{1+2\gamma}}\to_{L\to+\infty} +\infty \mbox{ when } \gamma>1/2,$$
$$\frac DL \ln\left(1+\frac M {D^2}\vee \sqrt{\frac M {D^2}}\right)\geq C(L/\ln L)^{\frac{-2\gamma}{1+2\gamma}},$$
and
$$\frac DM \geq  C (L/\ln L)^{\frac{-2\gamma}{1+2\gamma}}\mbox{ for }L \mbox{ large enough.}$$

{\bf Case 3.} If $\delta \leq \gamma  \leq 2\delta$, and $\delta \geq \gamma/(1+2\gamma)$, one chooses
$$M= 2^{J} \mbox{ with } J=\left\lfloor\log_2(L^{2/(1+4\delta)})\right \rfloor+1, $$ and
$$D=\left\lfloor M^{\delta/\gamma}\right \rfloor.$$

With such a choice, one has that $D\leq M$ and
\begin{eqnarray*}
R'^{2(1+2\gamma)}D^{-2\gamma}\geq C L^{\frac{-4\delta}{1+4\delta}},\\
R^2M^{-2\delta}\geq C L^{\frac{-4\delta}{1+4\delta}}.
\end{eqnarray*}
Furthermore
$$\frac{M}{D^2}\sim M^{1-2\delta/\gamma} \to_{L\to+\infty} 0 \mbox{ when } \delta>\gamma/2 \mbox{ and } 1 \mbox{ when } \delta=\gamma/2.$$
Hence,
$$\frac DL \ln\left(1+\frac M {D^2}\vee \sqrt{\frac M {D^2}}\right)\sim C\frac{\sqrt{M}}{L}\geq C L^{\frac{-4\delta}{1+4\delta}},$$
when $\delta\geq \gamma/2$, and
$$\frac DM \sim M^{\delta/\gamma -1}\geq L^{\frac{-4\delta}{1+4\delta}},$$
when $\delta\geq \gamma/(1+2\gamma)$.

{\bf Case 4.} If $\gamma\leq \delta$, one chooses  $M=D= 2^{J}$ with
$J=\left\lfloor\log_2( L^{2/(1+4\delta)})\right \rfloor+1$.

With such a choice,
$$R'^{2(1+2\gamma)}D^{-2\gamma}\geq R'^{2(1+2\gamma)}D^{-2\delta} \geq C L^{\frac{-4\delta}{1+4\delta}},$$
and
$$R^2M^{-2\delta}\geq C L^{\frac{-4\delta}{1+4\delta}}.$$

Moreover
$$\frac{M}{D^2} \to_{L\to+\infty} 0,$$
so
$$\frac DL \ln\left(1+\frac M {D^2}\vee \sqrt{\frac M {D^2}}\right)\sim \frac{\sqrt{M}}{L}\geq L^{\frac{-4\delta}{1+4\delta}},$$
and
$$\frac DM =1\geq L^{\frac{-4\delta}{1+4\delta}},$$
for $L$ large enough.

{\bf Case 5.} If $\gamma\leq 1/2$ and $\delta<\gamma/(1+2\gamma)$, one takes
$$M= 2^{J} \mbox{ with } J=\left\lfloor\log_2 L\right \rfloor+1, $$
and
$$D=\lfloor M^{1/(1+2\gamma)}\rfloor.$$

We first notice that $D\leq M$.
Then,
$$R'^{2(1+2\gamma)}D^{-2\gamma}\geq C L^{\frac{-2\gamma}{1+2\gamma}},$$
and
$$R^2M^{-2\delta}\geq CL^{-2\delta}.$$
Moreover
$$\frac{M}{D^2}\sim M^{\frac{2\gamma-1}{1+2\gamma}}\to_{L\to+\infty} 0  \mbox{ when } \gamma<1/2, \mbox{ and } 1 \mbox{ when } \gamma=1/2.$$
Hence,
$$\frac DL \ln\left(1+\frac M {D^2}\vee \sqrt{\frac M {D^2}}\right)\sim C\frac{\sqrt{M}}{L}\geq C L^{-1/2},$$
and
$$\frac DM \geq  C L^{\frac{-2\gamma}{1+2\gamma}}\geq C L^{-1/2}\mbox{ for }L \mbox{ large enough.}$$

This concludes the proof of Theorem \ref{minorations}.

\subsection{Proofs of Theorem \ref{puissancetest1} and Theorem \ref{puissancetest2}}\label{preuvetheor1et2}
\subsubsection{Preliminary results}

We consider here the general test function $\Phi_{\alpha}=\1_{\mathcal{T}_{\alpha}>0}$, defined by (\ref{fonctiontestgen}), where
$$\mathcal{T}_{\alpha}=\sup{\Lambda\in\mathcal{C}}\left(T_\Lambda''-t_{\Lambda,\alpha}''^{(N_L)}\right),$$
$T_\Lambda''=\sum_{\lambda\in \Lambda} T_\lambda$, and $\mathcal{C}$ is a finite collection of subsets of $\Lambda_\infty$. The collection $\mathcal{C}$ and the quantile  $t_{\Lambda,\alpha}''^{(N_L)}$ will be chosen to fit our two procedures respectively.

We begin to prove the following result.

\begin{Th}\label{puissancetestgene} Let $s\in\LL^{\infty}([0,1])$, and fix $\alpha$ and
$\beta$ in $]0,1[$. Assume that there exists some positive quantity $A_{\Lambda,\alpha,\beta}$ such that $$\PP_s\pa{t_{\Lambda,\alpha}''^{(N_L)}\geq A_{\Lambda,\alpha,\beta}}\leq \frac{\beta}{3}.$$
We recall that $D_\Lambda$ denotes the dimension of $S_\Lambda$ and we set $E_\Lambda=\sum_{j/(j,k)\in \Lambda} 2^j$.

There exist some positive constants
$C_1(\norm{s}_\infty,\beta)$ and $C_2(\beta)$ such that when $s$ satisfies

\begin{equation}
d^2(s,\calS_0) > \inf{\Lambda\in\mathcal{C}}\Bigg\{ \norm{s-s_\Lambda}^2
+C_1(\norm{s}_\infty,\beta)\pa{\frac{\sqrt{D_\Lambda}}{L}+\frac{\sqrt{E_\Lambda}}{L^{3/2}}}
+C_2(\beta)\frac{E_\Lambda}{L^2}+ A_{\Lambda,\alpha,\beta} \Bigg\}
\end{equation}
then
$$ \PP_s\pa{\Phi_{\alpha}=0}\leq \beta. $$

\end{Th}

\begin{proof}
Let $\alpha$ and $\beta$ in $]0,1[$, and $s$ be a fixed intensity.

\beqe \PP_s\pa{\Phi_{\alpha}=0}&=&\PP_s\pa{\mathcal{T}_{\alpha}\leq 0}\\
&=&\PP_s\pa{\forall \Lambda \in \mathcal{C}, T_\Lambda'' \leq t_{\Lambda,\alpha}''^{(N_L)}}\\
&\leq& \inf{\Lambda \in \mathcal{C}} \PP_s  \pa{ T_\Lambda'' \leq t_{\Lambda,\alpha}''^{(N_L)}}. \eeqe

For every $\Lambda$ in $\mathcal{C}$, we can write $T_\Lambda''$ in the following way : \beqe
T_\Lambda'' &=& \frac{1}{L^2} \sum_{\lambda \in \Lambda} \cro{\pa{\int_{[0,1]} \phi_\lambda(x)dN_x}^2-\int_{[0,1]}  \phi_\lambda^2(x)dN_x}\\
&=&  \frac{1}{L^2}  \sum_{\lambda \in \Lambda} \cro{\pa{\int_{[0,1]} \phi_\lambda(x) \pa{dN_x -s(x)L dx}}^2 +2 \int_{[0,1]}  \phi_\lambda(x) dN_x \int_{[0,1]}  \phi_\lambda(x) s(x) Ldx}\\
&& - \frac{1}{L^2}\sum_{\lambda \in \Lambda} \cro{\pa{\int_{[0,1]} \phi_\lambda(x) s(x) L dx}^2 + \int_{[0,1]}
\phi_\lambda^2(x)dN_x}. \eeqe

By setting
$$U_\Lambda=\frac{1}{L^2} \sum_{\lambda \in \Lambda} \cro{\pa{\int_{[0,1]} \phi_\lambda(x) \pa{dN_x -s(x) Ldx}}^2  - \int_{[0,1]}  \phi_\lambda^2(x)dN_x} $$
and
$$V_\Lambda= \frac{2}{L} \int_{[0,1]} (s_\Lambda(x)-\alpha_0\phi_0(x))\pa{dN_x -s(x) Ldx},$$
we obtain the following decomposition :
$$T_\Lambda'' = U_\Lambda+V_\Lambda+\|s_\Lambda\|^2-\alpha_0^2.$$
Since $d^2(s,\calS_0)=\|s-s_\Lambda\|^2+\|s_\Lambda\|^2-\alpha_0^2$,
 it follows that
$$ T_\Lambda''=U_\Lambda+V_\Lambda+d^2(s,\calS_0)-\|s-s_\Lambda\|^2.$$

Hence, \beq \label{decomposition} \PP_s\pa{\Phi_{\alpha}=0}\leq \inf{\Lambda\in \mathcal{C}} \PP_s  \pa{
U_\Lambda+V_\Lambda+d^2(s,\calS_0)\leq \|s-s_\Lambda\|^2 +t_{\Lambda,\alpha}''^{(N_L)}}.\eeq

The aim of the following lemmas is to define positive quantities $A_{\Lambda,\beta}^{(1)}$ and $A_{\Lambda,\beta}^{(2)}$,
 such that \beqe
&&\PP_s   \pa{ U_\Lambda\leq - A_{\Lambda,\beta}^{(1)} }\leq \frac{\beta}{3},\\
&&\PP_s   \pa{V_\Lambda \leq - A_{\Lambda,\beta}^{(2)} }\leq \frac{\beta}{3}.\\
\eeqe

Using (\ref{decomposition}) and assuming that $$\PP_s\pa{t_{\Lambda,\alpha}''^{(N_L)}\geq A_{\Lambda,\alpha,\beta}}\leq \frac{\beta}{3},$$ we then obtain  that as soon as there exists $\Lambda$ in $\mathcal{C}$ such that
\beq\label{distanceth}
 d^2(s,\calS_0) > \|s-s_\Lambda\|^2 + A_{\Lambda,\beta}^{(1)}+A_{\Lambda,\beta}^{(2)}+A_{\Lambda,\alpha,\beta} ,\eeq then
$$ \PP_s\pa{\Phi_{\alpha}=0}\leq \beta. $$

\begin{Lemma}\label{controlum} There exists some positive constant $C$ such that for all $\Lambda\in\mathcal{C}$ and for all $x>0$,
$$\PP_s\pa{-U_\Lambda\geq C\pa{\norm{s}_\infty \frac{\sqrt{D_\Lambda}}{L}\sqrt{x}+\norm{s}_\infty \frac{\sqrt{D_\Lambda}}{L} x+
\sqrt{\frac{\norm{s}_\infty E_\Lambda}{L^3}}x^{3/2}+\frac{E_\Lambda}{L^2}x^2}}\leq  2.77 e^{-x}.$$
\end{Lemma}

\begin{proof}
Let us first notice that \beqe U_\Lambda&=&\frac1{L^2}\sum_{\lambda\in \Lambda}\Bigg[\pa{\int_{[0,1]} \phi_\lambda(x)
dN_x}^2-\int_{[0,1]} \phi_\lambda^2(x) dN_x\\
&&-2 \pa{\int_{[0,1]}\phi_\lambda(x)dN_x}\pa{\int_{[0,1]} \phi_\lambda(x)s(x)Ldx}+ \pa{\int_{[0,1]}\phi_\lambda(x)s(x) Ldx}^2\Bigg]\\
&=&\frac1{L^2}\sum_{\lambda\in \Lambda}\Bigg[\sum_{l\neq l'=1}^{N_L}\phi_\lambda(X_l)\phi_\lambda(X_{l'})-2 \pa{\int_{[0,1]} \phi_\lambda(x)dN_x}\pa{\int_{[0,1]} \phi_\lambda(x)s(x) Ldx}\\
&&+ \pa{\int_{[0,1]}
\phi_\lambda(x)s(x) Ldx}^2\Bigg]\\
&=&\frac2{L^2}\sum_{\lambda\in \Lambda}\Bigg[\int_0^1\int_0^{y^-}
\phi_\lambda(x)\phi_\lambda(y) dN_x dN_y \\
&&-\pa{\int_0^1 \int_0^{y^-} \phi_\lambda(x)\phi_\lambda(y) dN_x s(y)Ldy + \int_0^1 \int_{y^{-}}^{1}
\phi_\lambda(x)\phi_\lambda(y) dN_x s(y)Ldy}\\
&&+\int_0^1 \int_0^{y^-}\phi_\lambda(x)\phi_\lambda(y) s(x) Ldx
s(y)Ldy\Bigg]\\
&=&\frac2{L^2}\sum_{\lambda\in \Lambda}\Bigg[\int_0^1\int_0^{y^-}
\phi_\lambda(x)\phi_\lambda(y) dN_x dN_y \\
&&-\pa{\int_0^1 \int_0^{y^-} \phi_\lambda(x)\phi_\lambda(y) dN_x s(y)Ldy + \int_0^1 \int_{0}^{x^{-}}
\phi_\lambda(x)\phi_\lambda(y) s(y)Ldy dN_x}+\\
&&\int_0^1 \int_0^{y^-}\phi_\lambda(x)\phi_\lambda(y) s(x) Ldx
s(y)Ldy\Bigg]\\
&=&\frac2{L^2}\sum_{\lambda\in \Lambda}\Bigg[\int_0^1\int_0^{y^-}\phi_\lambda(x)\phi_\lambda(y) (dN_x-s(x)Ldx)
(dN_y-s(y)Ldy)\Bigg].
\eeqe

Setting \beq \label{Hm} H_\Lambda(x,y)= \frac2{L^2}\sum_{\lambda\in \Lambda}\phi_\lambda(x)\phi_\lambda(y),\eeq we deduce
from Theorem 4.2 in Houdré and Reynaud-Bouret \cite{ustats} that there exists some absolute constant $\kappa>0$ such that for all $x>0$,
$$\PP_s\pa{-U_\Lambda\geq \kappa\pa{A_1\sqrt{x}+A_2 x+A_3x^{3/2}+A_4x^2}}\leq 2.77 e^{-x},$$ where

\beqe
A_1^2 &=&\int_0^1\int_0^y H_\Lambda^2(x,y)s(x)L dx s(y)Ldy,\\
A_2& =&\sup{a,b,\int a^2(x)s(x)Ldx=\int b^2(x)s(x)Ldx=1}\int_0^1 a(x) \pa{\int_x^1 b(y)H_\Lambda(x,y)s(y)Ldy} s(x)Ldx,\\
A_3^2&=& \sup{y\in[0,1]} \int_0^1 H_\Lambda^2(x,y) s(x)Ldx,\\
A_4 &=&\sup{x,y\in[0,1]} |H_\Lambda(x,y)|. \eeqe

Let us now evaluate $A_1$, $A_2$, $A_3$ and $A_4$ for every $\Lambda\in\mathcal{C}$.

To give an upper bound for $A_1^2$, we notice that $$A_1^2\leq \norm{s}_\infty^2 L^2 \int_0^1\int_0^1 H_\Lambda^2(x,y) dx
dy.$$

Since $\ac{\phi_\lambda,\lambda\in \Lambda}$ is an orthonormal basis on $[0,1]$, one has \beqe A_1^2&\leq
&\frac{4\norm{s}_\infty^2}{L^2} \int_0^1\int_0^1
 \sum_{\lambda\in \Lambda}\phi_\lambda(x)^2 \phi_\lambda(y)^2 dxdy\\
&\leq &\frac{4\norm{s}_\infty^2D_\Lambda}{L^2}. \eeqe

For $A_2$, we use Cauchy-Schwarz inequality to see that
$$A_2\leq \sup{b,\int b^2(x)s(x)Ldx=1}\cro{\int_0^1  \pa{\int_x^1 b(y)H_\Lambda(x,y)s(y)Ldy}^2 s(x) L dx}^{1/2},$$
and $$\pa{\int_x^1 b(y) H_\Lambda(x,y) s(y)Ldy}^2\leq \pa{\int_x^1 b^2(y)s^2(y) Ldy}\pa{\int_x^1 H_\Lambda^2(x,y) L dy}.$$
This implies that
\beqe A_2&\leq& L\cro{\int_0^1 \norm{s}_\infty\pa{\int_0^1 H_\Lambda^2(x,y)dy}s(x)dx}^{1/2}\\
&\leq &L\norm{s}_\infty \cro{\int_0^1 \int_0^1H_\Lambda^2(x,y) dxdy}^{1/2}\eeqe

Since $\ac{\phi_\lambda,\lambda\in \Lambda}$ is an orthonormal basis on $[0,1]$, one has

\beqe A_2 &\leq &\frac{2}{L}\norm{s}_\infty \cro{\int_0^1\int_0^1
 \sum_{\lambda\in \Lambda}\phi_\lambda^2(x)\phi_\lambda^2(y)dxdy}^{1/2}\\
&\leq & \frac{2}{L}\norm{s}_\infty D_\Lambda^{1/2} .\eeqe

As for $A_3$, we can prove in the same way that
$$A_3^2\leq\frac{4\norm{s}_\infty}{L^3}\sup{y\in[0,1]} \sum_{\lambda\in \Lambda} \phi_\lambda^2(y).$$
Moreover, for any $y$ fixed in $[0,1]$, \beqe
\sum_{(j,k)\in \Lambda} \phi_{j,k}^2(y) &\leq &
\sum_{j/(j,k)\in \Lambda} 2^{j}\\
&\leq & E_\Lambda.\eeqe

This implies that
$$A_3^2\leq \frac{4\norm{s}_\infty E_\Lambda}{L^3}.$$

Furthermore, for $x,y$ in [0,1],
\beqe |H_\Lambda(x,y)|&=&\frac{2}{L^2}\left| \sum_{(j,k)\in \Lambda}
\phi_{(j,k)}(x)
\phi_{(j,k)}(y)\right|\\
&\leq & \frac{2}{L^2}
\sum_{j/(j,k)\in \Lambda} 2^{j}\\
&\leq & \frac{2E_\Lambda}{L^2}.\eeqe

Finally, $A_4\leq {2E_\Lambda}/{L^2}$, and this concludes the proof of Lemma \ref{controlum}.
\end{proof}

By taking $x= \ln(8.31/\beta)$ in Lemma \ref{controlum}, we obtain that a possible value for $A_{\Lambda,\beta}^{(1)}$
is
$$A_{\Lambda,\beta}^{(1)}=C\pa{\norm{s}_\infty \frac{\sqrt{D_\Lambda}}{L} 2\ln(8.31/\beta)    +
\sqrt{\frac{\norm{s}_\infty E_\Lambda}{L^3}}(\ln(8.31/\beta))^{3/2}+\frac{E_\Lambda}{L^2}(\ln(8.31/\beta))^2},$$
where $C$ is an absolute positive constant.
We now use the following lemma, which derives from an analogue of Bennett's inequality (see proposition 7 of Reynaud-Bouret \cite{ptrfpois}, for instance).

\begin{Lemma}\label{controlvm}
There exists some positive constant $C$ such that for all $x>0$,
$$\PP_s\pa{-V_\Lambda\geq \frac{1}{2}\norm{s-\alpha_0\phi_0}^2-\frac{1}{2} \norm{s-s_\Lambda}^2 + \frac {C\norm{s}_\infty}{L} x}\leq e^{-x}.$$
\end{Lemma}

\begin{proof}
Recall that $$V_\Lambda= \frac{2}{L} \int_{[0,1]} (s_\Lambda(x)-\alpha_0\phi_0(x))\pa{dN_x -s(x) Ldx}.$$ Using proposition 7
of Reynaud-Bouret \cite{ptrfpois}, we easily obtain that for all $x>0$,
$$\PP\left(-V_\Lambda\geq   2
\sqrt{2x\frac{\norm{s}_\infty}{L}\norm{s_\Lambda-\alpha_0\phi_0}^2}+\frac{2\norm{s_\Lambda-\alpha_0\phi_0}_\infty}{3L}x\right)\leq
e^{-x}.$$

First note that
$$\norm{s_\Lambda-\alpha_0\phi_0}_\infty\leq \norm{s}_\infty.$$

By using the elementary inequality $2ab\leq {a^2}/{2}+2b^2$,
  we obtain that
\beqe
 2\sqrt{2x\frac{\norm{s}_\infty}{L}\norm{s_\Lambda-\alpha_0\phi_0}^2}&\leq &\frac{1}{2} \norm{s_\Lambda-\alpha_0\phi_0}^2 +
4 x\frac{\norm{s}_\infty}{L}\\
&\leq & \frac{1}{2}\norm{s-\alpha_0\phi_0}^2-\frac{1}{2}\norm{s-s_\Lambda}^2+4x\frac{\norm{s}_\infty}{L}. \eeqe

We deduce that there exists $C>0$ such that for all $x>0$,
$$\PP\left(-V_\Lambda\geq \frac{1}{2}\norm{s-\alpha_0\phi_0}^2 -\frac{1}{2} \norm{s-s_\Lambda}^2  + \frac {C\norm{s}_\infty}{L} x\right)\leq e^{-x}.$$
 \end{proof}

By taking $x= \ln(3/\beta)$ in Lemma \ref{controlvm}, we obtain that a possible value for $A_{\Lambda,\beta}^{(2)}$ is
$$A_{\Lambda,\beta}^{(2)}=\frac{1}{2}\norm{s-\alpha_0\phi_0}^2-\frac{1}{2} \norm{s-s_\Lambda}^2 + \frac {C\norm{s}_\infty}{L}\ln(3/\beta).$$

Replacing  $A_{\Lambda,\beta}^{(1)}$ and $A_{\Lambda,\beta}^{(2)}$ in (\ref{distanceth}) by the possible
values obtained above finally leads to the result of Theorem \ref{puissancetestgene}.
\end{proof}

\bigskip

We now prove the following lemma that will provide an upper bound for the quantity
$ A_{\Lambda,\alpha,\beta}$ occurring in Theorem \ref{puissancetestgene}.

\begin{Lemma}\label{controlquant}
Let $\tilde{X}_1, \ldots, \tilde{X}_n$ be  i.i.d. uniformly distributed on $[0,1]$. For $n\in \N $ and $\Lambda \subset \Lambda_{\infty}$, let
$$ T''_{\Lambda,n}=\frac1{L^2} \sum_{\lambda\in \Lambda} \sum_{l \neq l'=1}^n \phi_{\lambda}(\tilde{X}_l)\phi_{\lambda}( \tilde{X}_{l'}).$$
Let $D_\Lambda$ denote the dimension of $S_\Lambda$ and  $E_\Lambda=\sum_{j/(j,k)\in \Lambda} 2^j$.
 There exists some absolute constant $C>0$ such that for all $x>0$,
\begin{equation}\label{majquant}
\PP\pa{  T''_{\Lambda,n}\geq \frac{C n}{L^2} \pa{ \sqrt{D_{\Lambda}x}+x+\frac{E_{\Lambda}x^2}{n\vee 1}}}\leq 2.77 e^{-x}.
\end{equation}
 \end{Lemma}
\begin{proof}
If $n\in \{0,1\}$, $T''_{\Lambda,n}=0$ hence (\ref{majquant}) holds.
Since for all $\lambda \in \Lambda_{\infty}$, $\phi_{\lambda}$ is orthonormal to $\phi_0=\1_{[0,1]}$, it follows that the variables $ \phi_{\lambda}(\tilde{X}_l)$ are centered and we can  apply Theorem 3.4 in Houdré and Reynaud-Bouret \cite{ustats}. We now set $H_{\Lambda}(x,y)= \sum_{\lambda\in \Lambda}  \phi_{\lambda}(x)\phi_{\lambda}(y)/L^2.$ We obtain that there exists some absolute constant $C>0$ such that for all $x>0$,
$$\PP \pa{ T''_{\Lambda,n}  \geq C\pa{\tilde{A}_1\sqrt{x}+\tilde{A}_2 x+\tilde{A}_3x^{3/2}+\tilde{A}_4x^2}}\leq 2.77 e^{-x},$$ where

\beqe
\tilde{A}_1^2 &=&n^2\EE\cro{H_{\Lambda}^2(\tilde{X}_1,\tilde{X}_{2})}\\
\tilde{A}_2& =&\sup{}\ac{\left|\EE\cro{\sum_{l=1}^n \sum_{l'=1}^{l-1} H_{\Lambda}(\tilde{X}_1,\tilde{X}_{2})\alpha_l(\tilde{X}_1)\beta_{l'}(\tilde{X}_2)}\right|,\EE\cro{\sum_{l=1}^n \alpha_l^2(\tilde{X}_l)}\leq 1,\EE\cro{\sum_{l=1}^n \beta_l^2(\tilde{X}_l)}\leq 1}, \\
\tilde{A}_3^2&=& n\sup{y\in[0,1]} \int_0^1 H_{\Lambda}^2(x,y) dx,\\
\tilde{A}_4 &=&\sup{x,y\in[0,1]} |H_{\Lambda}(x,y)|. \eeqe

To evaluate $\tilde{A}_1$, $\tilde{A}_2$, $\tilde{A}_3$, $\tilde{A}_4$, we use arguments similar to the ones used in the
proof of Lemma \ref{controlum}.

Since $\ac{\phi_\lambda,\lambda\in \Lambda}$ is an orthonormal basis on $[0,1]$, \beqe \tilde{A}_1^2&\leq
& \frac{n^2}{L^4} \int_0^1\int_0^1
 \sum_{\lambda\in \Lambda}\phi_\lambda(x)^2 \phi_\lambda(y)^2 dxdy\\
&\leq &\frac{n^2D_{\Lambda}}{L^4}. \eeqe
Let $(\alpha_1,\ldots,\alpha_n)$ and $(\beta_1,\ldots,\beta_n)$ such that $\EE\cro{\sum_{l=1}^n
\alpha_l^2(\tilde{X}_l)}\leq 1$ and $\EE\cro{\sum_{l=1}^n \beta_l^2(\tilde{X}_l)}\leq 1$. Then

\beqe &&\left|\EE\cro{\sum_{l=1}^n \sum_{l'=1}^{l-1}
H_{\Lambda}(\tilde{X}_1,\tilde{X}_{2})\alpha_l(\tilde{X}_1)\beta_{l'}(\tilde{X}_2)}\right|\\
&&=\sum_{l=1}^n \sum_{l'=1}^{l-1}\int_0^1\int_0^1 H_{\Lambda}(x,y)\alpha_l(x)\beta_{l'}(y)dxdy\\
&&=\frac{1}{L^2}\sum_{l=1}^n \sum_{l'=1}^{l-1} \sum_{\lambda\in \Lambda} \int_0^1\phi_\lambda(x)\alpha_l(x) dx
\int_0^1\phi_\lambda(y)\beta_{l'}(y) dy.\eeqe

By using Cauchy-Schwarz inequality, we obtain \beqe &&\left|\EE\cro{\sum_{l=1}^n \sum_{l'=1}^{l-1}
H_{\Lambda}(\tilde{X}_1,\tilde{X}_{2})\alpha_l(\tilde{X}_1)\beta_{l'}(\tilde{X}_2)}\right|\\
&&\leq \frac{1}{L^2}\sum_{l=1}^n \sum_{l'=1}^{l-1} \cro{\sum_{\lambda\in \Lambda}
\pa{\int_0^1\phi_\lambda(x)\alpha_l(x) dx}^2}^{1/2}\cro{\sum_{\lambda\in \Lambda}
\pa{\int_0^1\phi_\lambda(y)\beta_{l'}(y) dy}^2}^{1/2}.\eeqe
One has for all $g\in\LL^2([0,1])$,
$\sum_{\lambda\in \Lambda} (\int \phi_\lambda g)^2\leq \int g^2$. As a consequence,
 \beqe
\tilde{A}_2&\leq &\frac{1}{L^2}\sum_{l=1}^n \cro{\int_0^1 \alpha_l^2(x) dx}^{1/2}\sum_{l'=1}^{n}\cro{\int_0^1
\beta_{l'}^2(y) dy}^{1/2}\\
&\leq&\frac{n}{L^2}.\eeqe

We evaluate $\tilde{A}_3^2$ and $\tilde{A}_4$ in the same way as $A_3^2$ and $A_4$ in the proof of Lemma
\ref{controlum}. We obtain that
$$\tilde{A}_3^2\leq \frac{ n E_{\Lambda}}{L^4},$$
and $$\tilde{A}_4\leq  E_{\Lambda} /L^2.$$
Finally, we proved that there exists some absolute constant $C>0$ such that
$$\PP\pa{ T''_{\Lambda,n} \geq C\frac{n}{L^2}
\pa{\sqrt{D_{\Lambda} x}+x+\frac{\sqrt{ E_{\Lambda} }}{\sqrt{n}}x^{3/2}+\frac{ E_{\Lambda}}{n}x^2}}\leq 2.77 e^{-x}.$$ Since
$$2\frac{\sqrt{E_{\Lambda}}}{\sqrt{n}}x^{3/2}\leq x+\frac{E_{\Lambda}}{n} x^2,$$
we can simplify the above inequality : there exists some constant $C>0$ such that
$$\PP\pa{  T''_{\Lambda,n} \geq C\frac{n}{L^2}
\pa{\sqrt{D_{\Lambda}  x}+x+\frac{E_{\Lambda} }{n}x^2}}\leq 2.77 e^{-x},$$
for all $x>0$. This concludes the proof of Lemma \ref{controlquant}.
\end{proof}
We are now in position to prove Theorem \ref{puissancetest1} and Theorem \ref{puissancetest2}.

\subsubsection{Proof of Theorem \ref{puissancetest1}}

Recall that the test function defined by (\ref{fonctiontest1}) is of the same form as the test function (\ref{fonctiontestgen}) of Theorem \ref{puissancetestgene} with $\mathcal{C}=\{\Lambda_J,\ J\in\mathcal{J}\}$, and $t_{\Lambda_J,\alpha}''^{(n)}=q_J'^{(n)}(u_{J,\alpha}'^{(n)})$, where $q_J'^{(n)}(u)$ denotes the $(1-u)$ quantile of $T''_{\Lambda_J,n}$.
 Since $u_{J,\alpha}'^{(n)}$ defined by (\ref{choixpoids}) satisfies  $u_{J,\alpha}'^{(n)}\geq \alpha e^{-W_J}$ for all $n$, one has that for all $n$,
 $$t_{\Lambda_J,\alpha}''^{(n)}\leq q_J'^{(n)}(\alpha e^{-W_J}).$$

In order to use Theorem \ref{puissancetestgene}, we then need to find some positive quantity $A_{J,\alpha,\beta}$ such that
\begin{equation}\label{AJalphabeta}
\PP_s\pa{q_J'^{(N_L)}(\alpha e^{-W_J})\geq A_{J,\alpha,\beta}}\leq \frac{\beta}{3}.
\end{equation}

Let us first give an upper bound for $q_J'^{(n)}(\alpha e^{-W_J})$ for all $n$ in $\N $.
We apply (\ref{majquant}) with $\Lambda=\Lambda_J$ (note that $ D_{\Lambda_J}=E_{\Lambda_J}=D_J$)
and with $x= \ln(2.77/\alpha)+W_J$. There exists some absolute constant $C>0$ such that
\begin{equation*}\label{evalq}
q_J'^{(n)}(\alpha e^{-W_J})\leq C \frac{n}{L^2}
\pa{\sqrt{D_J \pa{\ln(2.77/\alpha)+W_J}}+\ln(2.77/\alpha)+W_J+\frac{D_J}{n\vee 1}(\ln(2.77/\alpha)+W_J)^2}.
\end{equation*}

This allows  us to obtain some $A_{J,\alpha,\beta}$ such that (\ref{AJalphabeta}) holds. It actually gives that
$$q_J'^{(N_L)}(\alpha e^{-W_J})\leq C \frac{N_L}{L^2} \pa{\sqrt{D_J
\pa{\ln(2.77/\alpha)+W_J}}+\ln(2.77/\alpha)+W_J}+\frac{D_J}{L^2}(\ln(2.77/\alpha)+W_J)^2.$$
Now, from Bernstein's inequality, we deduce that for all $u>0$,
$$\PP_s \pa{N_L\geq \int_{[0,1]} s(x)Ldx + \sqrt{2\int_{[0,1]} s(x)Ldx u}+\frac{1}{3} u}\leq e^{-u}.$$
Hence a possible value for $A_{\Lambda_J,\alpha,\beta}$ is
$$C \frac{\int_{[0,1]} s(x)Ldx+\ln(3/\beta)}{L^2} \pa{\sqrt{D_J
\pa{\ln(2.77/\alpha)+W_J}}+\ln(2.77/\alpha)+W_J}+\frac{D_J}{L^2}(\ln(2.77/\alpha)+W_J)^2.$$

Using Theorem \ref{puissancetestgene} finally leads to the result of Theorem \ref{puissancetest1}.

\subsubsection{Proof of Theorem \ref{puissancetest2}}

Recall here that the test function defined by (\ref{fonctiontest2}) is of the same form as the test function (\ref{fonctiontestgen}) of Theorem \ref{puissancetestgene} with $\mathcal{C}=\{\Lambda,\ \Lambda\subset \Lambda_\bJ\}$, and $t_{\Lambda,\alpha}''^{(n)}=\sum_{\lambda\in\Lambda}q_\lambda^{(n)}\left({u_{\alpha}^{(n)}}/({2^j \bJ})\right)$, where $q_\lambda^{(n)}(u)$ denotes the $(1-u)$ quantile of $T_\lambda$ conditionally on the event $N_L=n$ under the
 null hypothesis $(H_0)$ and  $u_{\alpha}^{(n)}$ defined by (\ref{ualpha}) satisfies $u_{\alpha}^{(n)} \geq \alpha$ for all $n$.

Hence, we can prove Theorem \ref{puissancetest2} by using Theorem \ref{puissancetestgene} and some positive quantity $A_{\Lambda,\alpha,\beta}$ such that
 $$\PP_s\pa{\sum_{(j,k)\in \Lambda}q_{(j,k)}^{(N_L)}\left(\frac{\alpha}{2^{j} \bJ}\right)\geq A_{\Lambda,\alpha,\beta}}\leq \frac{\beta}{3}.$$
Following the same lines of proof as in the previous section, let us first give an upper bound for
 $q_{(j,k)}^{(n)}\left(\alpha/(2^j \bJ)\right)$.

Notice that $q_{(j,k)}^{(n)}(u)$ is  the
 $(1-u)$ quantile of the variable $T''_{\Gamma,n}$ with $\Gamma=\ac{(j,k)}$. Since $D_{\Gamma}=1$ and $E_{\Gamma}=2^j$, the inequality (\ref{majquant}) implies that there exists some constant $C>0$ such that for all $x>0$,
$$\PP \pa{ T''_{\Gamma,n} \geq  C\frac{n}{L^2} \pa{\sqrt{x}+x+2^j\frac{x^2}{n\vee 1}}}\leq 2.77 e^{-x}. $$
Taking $x=\ln(2.77)+\ln(2^j\bJ/\alpha)$ in this inequality leads to the conclusion that :
$$q_{(j,k)}^{(n)}\left(\frac{\alpha}{2^j\bJ}\right)\leq C\frac{n}{L^2}
\pa{\sqrt{\pa{\ln(2.77)+\ln(2^j\bJ/\alpha)}}+\ln(2.77)+\ln(2^j\bJ/\alpha)}+\frac{2^j}{L^2}(\ln(2.77)+\ln(2^j\bJ/\alpha))^2.$$

From Bernstein's inequality, we deduce that  a possible value
for $A_{\Lambda,\alpha,\beta}$ is

\begin{multline}
C\sum_{(j,k)\in \Lambda} \Bigg\{\frac{\int_{[0,1]} s(x)Ldx+\ln(3/\beta)}{L^2}
\pa{\sqrt{\pa{\ln(2.77)+\ln(2^j\bJ/\alpha)}}+\ln(2.77)+\ln(2^j\bJ/\alpha)}\\
+\frac{2^j}{L^2}(\ln(2.77)+\ln(2^j\bJ/\alpha))^2\Bigg\},
\end{multline}
for some positive constant $C$.

Since $ |\Lambda|=D_\Lambda-1$ and $E_{\Lambda}\leq 2^\bJ$, we obtain the result
of Theorem \ref{puissancetest2}.\\

\subsection{Proof of Proposition \ref{vitesses1}}

Let us assume that $s$ belongs to $\mathcal{B}_{2,\infty}^\delta(R)\cap  \LL^\infty(R'')$. We need to find an upper bound for
the quantity
\begin{multline*}
\inf{J \in \mathcal{J}} \Bigg\{
\norm{s-s_J}^2
+C_1(\norm{s}_\infty,\beta)\frac{\sqrt{D_J}}{L}+C_2(\beta)\frac{D_J}{L^2}
+C_3(\alpha,\beta)\int_{[0,1]}
s(x)dx\pa{\frac{\sqrt{D_JW_J}}{L}+\frac{W_J}{L}} \\+ C_4(\alpha)\frac{D_J
W_J^2}{L^2}\Bigg\},
\end{multline*}
in Theorem \ref{puissancetest1}.

We have already noticed that when $s$ belongs to $\mathcal{B}_{2,\infty}^\delta(R)$, for all $J\geq 1$, $$\|s-s_J\|^2 \leq c(\delta) R^2 D_J^{-2\delta}.$$ Moreover, the constant
$C_1(\|s\|_\infty,\beta)$ can be replaced by  $C_1(R'',\beta)$, so we only need to find an upper bound for
$$C(\alpha,\beta,R'',\delta)\inf{J \in \mathcal{J}} \Bigg\{
 R^2
D_J^{-2\delta}+\frac{\sqrt{D_J}}{L}+\frac{D_J}{L^2}+\frac{\sqrt{D_JW_J}}{L}+\frac{W_J}{L}
+ \frac{D_JW_J^2}{L^2} \Bigg\}.$$

Taking $W_J=\ln |\mathcal{J}|=\ln \left\lfloor\log_2(L^2/(\ln\ln L)^3)\right \rfloor$, with $\ln \ln L\geq 1$ leads to $W_J\leq 2.06 \ln\ln L$, so
\begin{multline*}
C(\alpha,\beta,R'',\delta)\inf{J \in \mathcal{J}} \Bigg\{
 R^2
D_J^{-2\delta}+\frac{\sqrt{D_J}}{L}+\frac{D_J}{L^2}+\frac{\sqrt{D_JW_J}}{L}+\frac{W_J}{L}
+ \frac{D_JW_J^2}{L^2} \Bigg\}\\
\leq C'(\alpha,\beta,R'',\delta)\left(\inf{J \in \mathcal{J}} \left\{
 R^2
D_J^{-2\delta}+\frac{\sqrt{D_J\ln\ln L}}{L}+\frac{D_J(\ln\ln L)^2}{L^2}\right\}+\frac{\ln\ln L}{L}\right).
\end{multline*}
Since for all $J$ in $\mathcal{J}$, $D_J\leq {L^2}/{(\ln\ln L)^3}$,
\begin{multline*}
C'(\alpha,\beta,R'',\delta)\left(\inf{J \in \mathcal{J}} \left\{
 R^2
D_J^{-2\delta}+\frac{\sqrt{D_J\ln\ln L}}{L}+\frac{D_J(\ln\ln L)^2}{L^2}\right\}+\frac{\ln\ln L}{L}\right)\\
\leq C''(\alpha,\beta,R'',\delta)\left(\inf{J \in \mathcal{J}} \left\{
 R^2
D_J^{-2\delta}+\frac{\sqrt{D_J\ln\ln L}}{L}\right\}+\frac{\ln\ln L}{L}\right).
\end{multline*}

We have that $R^2D_J^{-2\delta}<{\sqrt{D_J \ln\ln L}}/{L}$ if and
only if $J>\log_2\left(\left({R^4 L^2}/{\ln \ln L}\right)^{{1}/\pa{1+4\delta}}\right)$. Hence, we
introduce
 $$J_*=\left\lfloor\log_2 \left( \left( \frac{R^4 L^2}{\ln \ln
          L}\right)^{\frac{1}{1+4\delta}}\right)\right \rfloor +1,$$ and we distinguish three cases.\\
When $1\leq J_*\leq
          \left\lfloor\log_2(L^2/(\ln \ln L)^3)\right \rfloor,$ then $J_*$ belongs to $\mathcal{J}$ and
\begin{eqnarray*}
\inf{J \in \mathcal{J}} \Bigg\{R^2
D_J^{-2\delta}+\frac{\sqrt{D_J \ln\ln L}}{L}\Bigg\}&\leq& R^2
D_{J_*}^{-2\delta}+\frac{\sqrt{D_{J_*} \ln\ln L}}{L}\\
&\leq & (1+\sqrt{2})R^{\frac{2}{4\delta+1}} \left( \frac{\sqrt{\ln \ln L}}{L}
  \right)^{\frac{4\delta}{4\delta+1}}.
\end{eqnarray*}

When $J_*>\left\lfloor\log_2(L^2/(\ln \ln L)^3)\right \rfloor$, this means that for all $J$ in $\mathcal{J}$,
$\sqrt{D_J\ln \ln L}/L\leq R^2 D_J^{-2\delta}.$
By taking $J^*=\left\lfloor\log_2(L^2/(\ln \ln L)^3)\right \rfloor$, we obtain that
$$\inf{J \in \mathcal{J}} \Bigg\{R^2
D_J^{-2\delta}+\frac{\sqrt{D_J \ln\ln L}}{L}\Bigg\}\leq2R^2
D_{J^*}^{-2\delta}\leq 2^{2\delta+1}R^2\left(\frac{(\ln \ln L)^3}{L^2}\right)^{2\delta}.$$
Finally, when $J_*<1$, then for all $J$ in $\mathcal{J}$, $R^2D_J^{-2\delta}\leq \sqrt{D_J \ln\ln
    L}/L$, so by taking $J^*=1$, we obtain that
$$
\inf{J \in \mathcal{J}} \Bigg\{R^2
D_J^{-2\delta}+\frac{\sqrt{D_J \ln\ln L}}{L}\Bigg\}\leq 2\frac{\sqrt{2\ln\ln
  L}}{L}.$$
This ends the proof.

\subsection{Proof of Proposition \ref{vitesses2}}

Let us assume that $s$ belongs to $\mathcal{B}_{2,\infty}^\delta(R)\cap W_\gamma(R')\cap \LL^\infty(R'')$. We now need to find an adequate upper bound for
\begin{multline*}
\inf{\Lambda\subset\Lambda_\bJ} \Bigg\{ \norm{s-s_\Lambda}^2
+C_1(\norm{s}_\infty,\beta)\pa{\frac{\sqrt{D_\Lambda}}{L}+\frac{2^{\bJ/2}}{L^{3/2}}}+C_2(\beta)\frac{2^\bJ}{L^2}\\
+C_3(\alpha,\beta)\int_{[0,1]} s(x)dx\frac{D_\Lambda\ln(2^{\bJ}\bJ)}{L} + C_4(\alpha)\frac{D_\Lambda 2^\bJ
\ln^2(2^\bJ \bJ)}{L^2}\Bigg\}
\end{multline*}
in Theorem \ref{puissancetest2}.

As in the proof of Proposition \ref{vitesses1}, the constant
$C_1(\|s\|_\infty,\beta)$ can be replaced by  $C_1(R'',\beta)$. Moreover, with the choice $\bJ=\left\lfloor\log_2(L/\ln L)\right \rfloor$, we have that
$$\frac{2^{\bJ/2}}{L^{3/2}}\leq \frac{1}{L\sqrt{\ln L}},$$
$$\frac{2^\bJ}{L^2}\leq  \frac{1}{L\ln L},$$
and
 $\ln(2^\bJ\bJ)\leq \ln L$. So we only need to find an upper bound for
$$C(\alpha,\beta,R'')\left(\inf{\Lambda\in \Lambda_\bJ} \Bigg\{\norm{s-s_\Lambda}^2
+\frac{D_\Lambda\ln L}{L}\Bigg\}+\frac{1}{L\sqrt{\ln L}}\right).$$

Let us introduce for all integer $D\leq 2^\bJ$ the subset $\tilde{\Lambda}_D$ of $\Lambda_\bJ$ such that the elements of $\{\alpha_\lambda,\lambda\in \tilde{\Lambda}_D\}$  are the $(D-1)$ largest elements in $\{\alpha_\lambda, \lambda\in\Lambda_\bJ\}$.

We can notice that
$$\norm{s-s_{\tilde{\Lambda}_D}}^2=\norm{s-s_{\bJ}}^2+\norm{s_{\bJ}-s_{\tilde{\Lambda}_D}}^2.$$

On the one hand, since $s$ belongs to $\mathcal{B}_{2,\infty}^\delta(R)$, $$\|s-s_\bJ\|^2 \leq C(\delta) R^2 \left(\frac{L}{\ln L}\right)^{-2\delta}.$$

On the other hand, since $s$ belongs to $W_\gamma(R')$, then for all $t>0$,
\begin{eqnarray*}
\sum_{j\in\N} \sum_{k=0}^{2^j-1} \1_{|\alpha_{(j,k)}|>\frac{t}{2}}&\leq & \sum_{j\in\N} \sum_{k=0}^{2^j-1} \sum_{l\in\N} \1_{\frac{t}{2} 2^{l}< |\alpha_{(j,k)}|\leq \frac{t}{2} 2^{l+1}}\\
&\leq &\sum_{l\in\N} \sum_{j\in\N} \sum_{k=0}^{2^j-1}  \left(\frac{|\alpha_{(j,k)}|}{\frac{t}{2} 2^{l}}\right)^2 \1_{|\alpha_{(j,k)}|\leq \frac{t}{2} 2^{l+1}}\\
&\leq & 4\sum_{l\in\N} \frac{2^{-2l}}{t^2}\sum_{j\in\N} \sum_{k=0}^{2^j-1}  \alpha_{(j,k)}^2 \1_{|\alpha_{(j,k)}|\leq t 2^{l}}\\
&\leq & 4\sum_{l\in\N} \frac{2^{-2l}}{t^2}R'^2\left(t^22^{2l}\right)^{\frac{2\gamma}{1+2\gamma}}\\
&\leq & C(\gamma)R'^2 t^{-\frac{2}{1+2\gamma}}.
\end{eqnarray*}

Taking $t$ such that $C(\gamma)R'^2 t^{-\frac{2}{1+2\gamma}}=D$ in the above inequality proves that all the coefficients of $s_\bJ-s_{\tilde{\Lambda}_D}$ are smaller than $t/2$ and
\begin{eqnarray*}
\norm{s_\bJ-s_{\tilde{\Lambda}_D}}^2&\leq& \sum_{j=0}^{\bJ-1}\sum_{k=0}^{2^j-1} \alpha_{(j,k)}^2\1_{|\alpha_{(j,k)}|\leq \frac{t}{2}},\\
&\leq & C(\gamma)R'^{2+4\gamma} D^{-2\gamma}.
\end{eqnarray*}
Hence,
\begin{multline*}
C(\alpha,\beta,R'')\left(\inf{\Lambda\in \Lambda_\bJ} \Bigg\{\norm{s-s_\Lambda}^2
+\frac{D_\Lambda\ln L}{L}\Bigg\}+\frac{1}{L\sqrt{\ln L}}\right)\\
\leq C(\alpha,\beta,R'',\delta,\gamma)\left(\inf{1\leq D\leq 2^\bJ} \Bigg\{ R'^{2+4\gamma} D^{-2\gamma}+\frac{D\ln L}{L}\Bigg\}+R^2 \left(\frac{L}{\ln L}\right)^{-2\delta}+\frac{1}{L\sqrt{\ln L}}\right).
\end{multline*}

We have that $R'^{2+4\gamma}D^{-2\gamma}<D \ln L/L$ if and
only if $D>  R'^2\left( L/\ln L\right)^{1/\pa{1+2\gamma}}$. Hence, we
introduce
 $$D_*=\left\lfloor R'^2\left( \frac{L}{\ln L}\right)^{\frac{1}{1+2\gamma}}\right \rfloor+1,$$ and we distinguish two cases.\\

When $1\leq D_*\leq 2^\bJ$, we clearly obtain that
$$\inf{1\leq D\leq 2^\bJ} \Bigg\{ R'^{2+4\gamma} D^{-2\gamma}+\frac{D\ln L}{L}\Bigg\}\leq R'^{2+4\gamma} D_*^{-2\gamma}+\frac{D_*\ln L}{L}.$$

On the one hand, when $D_*\geq 2$, this leads to

\begin{eqnarray*}
\inf{1\leq D\leq 2^\bJ} \Bigg\{ R'^{2+4\gamma} D^{-2\gamma}+\frac{D\ln L}{L}\Bigg\}&\leq &R'^{2+4\gamma} D_*^{-2\gamma}+2\frac{(D_*-1)\ln L}{L}\\
&\leq& 3R'^2\left(\frac{L}{\ln L}\right)^{\frac{-2\gamma}{1+2\gamma}}.
\end{eqnarray*}

On the other hand, when $D_*=1$, since $R'^{2+4\gamma}D_*^{-2\gamma}<D_* \ln L/L$, one has

$$\inf{1\leq D\leq 2^\bJ} \Bigg\{ R'^{2+4\gamma} D^{-2\gamma}+\frac{D\ln L}{L}\Bigg\}\leq 2\frac{\ln L}{L}.$$

Now,  let us consider the case where $D_*>2^\bJ$. This means that for all $D$ such that $1\leq D\leq 2^\bJ$,
$D\ln L/L\leq R'^{2+4\gamma} D^{-2\gamma}.$
By taking $D^*=2^\bJ$, we obtain that
\begin{eqnarray*}
\inf{1\leq D\leq 2^\bJ} \Bigg\{ R'^{2+4\gamma} D^{-2\gamma}+\frac{D\ln L}{L}\Bigg\}&\leq &2R'^{2+4\gamma} {D^*}^{-2\gamma}\\
&\leq & 2R'^{2+4\gamma} \left(\frac{L}{2\ln L}\right)^{-2\gamma}.
\end{eqnarray*}
This concludes the proof of Proposition \ref{vitesses2}.\\
\ \\
{\bf Acknowledgment.} The authors acknowledge the support of the French Agence Nationale de la Recherche (ANR), under grant ATLAS (JCJC06$\_$137446) ''From Applications to Theory in Learning and Adaptive Statistics''.

\end{document}